\documentclass[notitlepage,11pt]{amsart}
\usepackage{amssymb}
\usepackage{amsmath}
\usepackage{amscd}
\usepackage{xspace}
\usepackage[mathscr]{eucal}

\DeclareMathOperator{\hhh}{\ensuremath{Hom}\xspace}
\DeclareMathOperator{\eee}{\ensuremath{Ext}\xspace}
\DeclareMathOperator{\image}{\ensuremath{im}\xspace}
\DeclareMathOperator{\id}{\ensuremath{id}\xspace}

\DeclareMathOperator{\sn}{\ensuremath{sn}\xspace}
\DeclareMathOperator{\Rad}{\ensuremath{rad}\xspace}

\newcommand{\mapdef}[1]{\ensuremath{\overset{#1}{\longrightarrow}}\xspace}
\newcommand{\mf}[1]{\ensuremath{\mathfrak{#1}}\xspace}
\newcommand{\spn}{\ensuremath{\mf{sp}(2n)}\xspace}
\newcommand{\spl}{\ensuremath{\mf{sl}_2}\xspace}
\newcommand{\uu}{\ensuremath{\mf{U}(\spn)}\xspace}
\newcommand{\usl}{\ensuremath{\mf{U}(\spl)}\xspace}

\newcommand{\dd}{\ensuremath{\Delta}\xspace}
\newcommand{\dn}{\ensuremath{\Delta_0}\xspace}
\newcommand{\la}{\ensuremath{\lambda}\xspace}
\newcommand{\La}{\ensuremath{\Lambda}\xspace}
\newcommand{\aaa}{\ensuremath{\alpha}\xspace}
\newcommand{\nn}{\ensuremath{\mathbb{N}_0}\xspace}
\newcommand{\N}{\ensuremath{\mathbb N}\xspace}
\newcommand{\Z}{\ensuremath{\mathbb Z}\xspace}
\newcommand{\C}{\ensuremath{\mathbb C}\xspace}
\newcommand{\halfz}{\ensuremath{{\mathbb Z}_2}\xspace}
\newcommand{\calc}{\ensuremath{\mathcal{C}}\xspace}
\newcommand{\calh}{\ensuremath{\mathcal{H}}\xspace}
\newcommand{\calo}{\ensuremath{\mathcal{O}}\xspace}
\newcommand{\caloo}{\ensuremath{\mathcal{O}_\N}\xspace}
\newcommand{\calp}{\ensuremath{\mathcal{P}}\xspace}
\newcommand{\cals}{\ensuremath{\mathcal{S}}\xspace}
\newcommand{\cala}{\ensuremath{\mathcal{A}_0}\xspace}
\newcommand{\pfilt}{\ensuremath{\mathcal{F}(\dd)}\xspace}
\newcommand{\qfilt}{\ensuremath{\mathcal{F}(\nabla)}\xspace}
\newcommand{\Hom}{\ensuremath{\hhh_{A}}\xspace}
\newcommand{\HHom}{\ensuremath{\hhh_{H_f}}\xspace}
\newcommand{\hhom}{\ensuremath{\hhh_\calo}\xspace}
\newcommand{\Ext}{\ensuremath{\eee_\calo}\xspace}

\newcommand{\half}{\ensuremath{\frac{1}{2}}\xspace}

\newcommand{\start}[2]{\ensuremath{\varphi \big(W' \big( #1,#2 \big)
\big)\ :\ }\xspace}

\newcommand{\diag}[4]{\ensuremath \bigg[{#1
\overset{\sigma_1}{\longleftrightarrow} #2
\overset{\sigma_2}{\longleftrightarrow} #3
\overset{\sigma_1}{\longleftrightarrow} #4
\overset{\sigma_2}{\longleftrightarrow} #1} \bigg]\xspace}

\theoremstyle{plain}
\newtheorem{theorem}{Theorem}
\newtheorem{prop}{Proposition}
\newtheorem{lemma}{Lemma}
\newtheorem{cor}{Corollary}

\theoremstyle{definition}
\newtheorem{remark}{Remark}
\newtheorem{stand}{Standing Assumption}

\begin{document}
\title[Category \calo over a deformation of the symplectic oscillator
algebra]
{Category \calo over a deformation of the\\ symplectic oscillator
algebra}
\author[Apoorva Khare]{}
\date{\today}
\renewcommand{\thefootnote}{}
\footnote{{\it Email address}: {\tt apoorva@math.uchicago.edu} (Apoorva
Khare).}
\maketitle

\begin{center}
{\large Apoorva Khare}\\
\vspace{1ex}
{\small \it Department of Mathematics (University of Chicago),\\ 5734 S.
University Avenue, Chicago, IL - 60637, USA}
\end{center}

\begin{abstract}
\small{We discuss the representation theory of $H_f$, which is a
deformation of the symplectic oscillator algebra $\mf{sp}(2n) \ltimes
\mf{h}_n$, where $\mf{h}_n$ is the ($(2n+1)$-dimensional) Heisenberg
algebra. We first look at a more general algebra with a triangular
decomposition. Assuming the PBW theorem, and one other hypothesis, we
show that the BGG category \calo is abelian, finite length, and
self-dual.

We decompose \calo as a direct sum of blocks $\calo(\la)$, and show that
each block is a highest weight category.

In the second part, we focus on the case $H_f$ for $n=1$, where we prove
all these assumptions, as well as the PBW theorem.}
\end{abstract}

\tableofcontents

\section*{Introduction}

We discuss here the BGG category over a deformation of a well known
algebra $H_0 = \mf{U}(\mf{sp}(2n)) \ltimes A_n$. The relation $[Y_i,X_i]
= 1$ in $A_n$ is deformed using the quadratic Casimir operator $\dd$ of
$\mf{sp}(2n)$. We work throughout over a ground field $k$ of
characteristic zero.\\

In the first half, we work in a more general setup, involving an algebra
with a triangular decomposition. We carry out many of the classical
constructions, including standard (Verma) and co-standard modules, and
introduce the BGG category \calo. Next, we introduce the duality functor,
which is exact, and show some homological properties of \calo. Assuming
the non-vanishing and finite length of all Verma modules, we show that
\calo has many good properties (in particular, it is abelian, finite
length, and self-dual).

Under additional assumptions, we decompose \calo as a direct sum of
subcategories  - or {\it blocks} - $\calo(\la)$. We show that each of
these blocks $\calo(\la)$ - and hence $\calo$ - has enough projectives.
This helps us construct projective covers, injective hulls, and
progenerators in each block. There is also an equivalence from
$\calo(\la)$ to the category of finitely generated modules over a
finite-dimensional algebra. Assuming the PBW theorem, each block is a
highest weight category, so that BGG reciprocity holds here.\\

In the second half, we introduce our algebra $H_0$ (and later on, $H_f$),
and produce explicit automorphisms and an anti-involution (which is used
to consider duality). We then focus on the case $n=1$. Analogous to
\spl-theory, we first look at standard cyclic modules via explicit
calculations. We then show that a large set of Verma modules are
nonzero.

Next, we show that an important constant $\aaa_{rm}$ is actually a
polynomial. This shows the PBW Theorem. We then take a closer look at
Verma modules. There is an important condition for a Verma module $Z(r)$
to have a submodule $Z(t)$: the constant $\aaa_{r,r-t+1}$ above must
vanish. This helps partition $k$ into the blocks $S(r)$.

The structure of finite-dimensional simple modules is very similar to the
\spl-case; we state the well-known character formulae here. We completely
classify all Vermas with non-integer weights, and give some results on
Vermas with integer weights. Therefore all the assumptions (and results)
of the first half are shown to hold for $H_f$.\\

\section*{\bf Part 1 : General theory}

In this first part, we examine in detail the structure of the category
\calo, and several duality and homological properties, under a general
setup involving a general algebra with a triangular decomposition. (In
particular, this treatment is valid for a finite-dimensional semisimple
Lie algebra $\mf{g}$ over $\C$.) We end by showing that the category is a
direct sum of blocks, each of which is a (``finite-dimensional") highest
weight category. The main goal of the second part, will be to prove (for
the algebra $H_f$) the assumptions used in this part (including the PBW
theorem), so that the results proved here all hold. Thus, one may read
the second part independently of the first.\\

\section{Standard cyclic modules in the Harish-Chandra (or BGG)
category}\label{sec5}

\noindent{\bf Setup :}

\noindent We work throughout over a ground field $k$ of characteristic
zero. We define $\nn = \N \cup \{0\}$. We work over an associative
$k$-algebra $A$, having the following properties.

\begin{enumerate}
\item The multiplication map $: B_- \otimes_k H \otimes_k B_+
\twoheadrightarrow A$ is surjective, where all symbols denote associative
$k$-subalgebras of $A$ (this is the {\it triangular decomposition}).\\

\item There is a finite-dimensional subspace $\mf{h}$ of $H$ so that $H =
\text{Sym}(\mf{h})$. Thus $\mf{h}$ is an abelian Lie algebra (or $H$ is
abelian).\\

\item There exists a {\it base of simple roots} $\dd$, i.e. a basis $\dd$
of $\mf{h}^* = \hhh_k(\mf{h},k)$. Define a partial ordering on $\mf{h}^*$
by: $\la \geq \mu$ iff $\la - \mu \in \nn \dd$, i.e. $\la - \mu$ is a sum
of finitely many elements of $\dd$ (repetitions allowed).\\

\item $A = \bigoplus_{\mu \in \Z \dd} A_\mu$, where $A_\mu$ is a weight
space for ad $\mf{h}$. In other words, $[h,a_\mu] = h a_\mu - a_\mu h =
\mu(h) a_\mu$ for all $h \in \mf{h},\ \mu \in \Z \dd,\ a_\mu \in
A_\mu$. Further, $B_+ \subset \bigoplus_{\mu \in \nn \dd} A_\mu$ and $H
\subset A_0$.\\

\item $(B_+)_0 = k$, and $\dim_k(B_+)_\mu < \infty$ for every $\mu$.\\

\item There exists an anti-involution $i$ of $A$ (i.e. $i^2|_A = \id|_A$)
that takes $(B_+)_\mu$ to $(B_-)_{-\mu}$ for each $\mu$, and acts as the
identity on all of $H$.\\
\end{enumerate}

\begin{remark}
Because of the anti-involution $i$, similar properties are true for
$B_-$, as are mentioned for $B_+$ above. We also have subalgebras
(actually, ideals) $N_+ = \bigoplus_{\mu \neq 0} (B_+)_\mu$ in $B_+$, and
similarly, $N_-$ in $B_-$.\\
\end{remark}

For an ($A$- or) $H$-module $V$, denote by $\Pi(V)$ the set of {\it
weights} $\mu \in \mf{h}^*$, so that the {\it weight space} $V_\mu := \{
v \in V : h v = \mu(h) v\ \forall h \in \mf{h} \}$ is nonzero. Then
standard arguments say that $\sum_{\mu \in \mf{h}^*} V_\mu =
\bigoplus_{\mu \in \mf{h}^*} V_\mu$ is the largest $\mf{h}$-semisimple
submodule of $V$.

We now introduce the {\it Harish-Chandra category} \calh. Its objects are
$A$-modules with a (simultaneously) diagonalizable $\mf{h}$-action, and
finite-dimensional weight spaces. Clearly, \calh is a full abelian
subcategory of $A$-mod. Inside this, we also introduce the (full) BGG
subcategory \calo, whose objects are finitely generated objects of \calh
with a locally finite action of $B_+$, i.e. $\forall M \in \calo,\ B_+ m$
is finite-dimensional for each $m \in M$. Note that $\mathcal{O}$ is not
extension-closed in $A$-mod (cf. \cite{K}).\\

\noindent{\bf Definitions :}

A {\it maximal vector} in an $A$-module $V$ is a weight vector for
$\mf{h}$ that is killed by $N_+$; in other words, it is an eigenvector
for $B_+$.

A {\it standard cyclic module} is an $A$-module generated by exactly one
maximal vector. Certain universal standard cyclic modules are called {\it
Verma modules}, just as in the classical case of \cite{BGG1} or
\cite{H}.\\

There exist maximal vectors (i.e. eigenvectors for $B_+$) in any object
of \calo. We now look at standard cyclic modules, namely $V = A \cdot
v_\la$, where $v_\la$ is maximal with weight $\la$. Most (if not all) of
the results in \cite[$\S 20 \S$]{H} now hold. We can construct standard
cyclic modules $B_- v_\la$ and Verma modules $Z(\la) = A / (N_+,\{(h -
\la(h) \cdot 1) : h \in \mf{h} \})$ with unique simple quotients
$V(\la)$, for each $\la \in \mf{h}^*$.

\begin{stand}
Until Section \ref{hfn1}, we keep the assumption that every Verma module
$Z(\la)$ is nonzero.\\
\end{stand}

The $V(\la)$'s are pairwise non-isomorphic, exhaust all simple objects in
\calo, and are in bijective correspondence with $\mf{h}^*$, as well as
each of the sets of finite-dimensional simple $\mf{h}$-modules, and
finite-dimensional simple $(H \otimes_k B_+)$-modules. (For the last two
bijective correspondences, we also need $k$ to be algebraically closed,
so that we can use Lie's theorem. For the same reason, all
finite-dimensional simple modules are also in \calo, whenever $k$ is also
algebraically closed.)\\

\noindent {\it Notation:} Any standard cyclic module $V$ of highest
weight $\la$ is a quotient of $Z(\la)$. We denote this (or $V$) by
$Z(\la) \to V \to 0$. We also denote the annihilator of $V(\la)$ in $A$
by $J(\la)$, and the (unique) maximal submodule of $Z(\la)$ by $Y(\la)$,
so that $V(\la) = A / J(\la) = Z(\la) / Y(\la)$.

\begin{theorem}\label{2tfae}
Suppose $V \in \calo$. Then the following are equivalent:
\begin{enumerate}
\item $\Hom(Z(\la),V) \neq 0$.
\item $V$ has a maximal weight vector $v_\la$ of weight $\la$.
\item $V$ has a standard cyclic submodule $V'$ of highest weight $\la$.\\
\end{enumerate}
\end{theorem}

\noindent Now, by seeing where the maximal vector goes, we also have

\begin{cor}\label{C2.1}
If $Z(\la) \to V \to 0$, then $\dim_k (\Hom(V,V(\mu))) = \delta_{\mu \la}
\in \{ 0,1 \}$.
\end{cor}

\begin{lemma}\label{Ltfae}
If $V$ and $V'$ are standard cyclic of highest weight $\la$, then the
following are equivalent:
\begin{enumerate}
\item $V \to V' \to 0$.
\item $\Hom(V,V') = k$.
\item $\Hom(V,V') \neq 0$.\\
\end{enumerate}
\end{lemma}

We now define the {\it formal character} (cf. \cite[$\S 13,21 \S$]{H}) of
an $A$-module $V = \bigoplus_\mu V_\mu \in \calh$. This is just the
formal sum $ch_V = \sum_{\mu \in \mf{h}^*} (\dim V_\mu) e(\mu)$, where
$\Z[\mf{h}^*] = \bigoplus_{\la \in \mf{h}^*} \Z \cdot e(\la)$. Finally,
define the {\it Kostant function} $p(\la)$ to be $p(\la) =
\dim_k(B_+)_{-\la} = \dim_k(B_-)_\la$.\\

\section{Duality and homological properties}

As is standard, we give $M^* = \hhh_k(M,k)$ a {\it left} $A$-module
structure (for each $M \in \calo$), using the anti-involution $i$
mentioned above. Now define the functor $F$ from $\calo$ to the opposite
category $\calo^{op}$ (defined presently), by taking $F(M)$ to be the
submodule of $M^*$ generated by all $\mf{h}$-weight vectors in $M^*$.
Thus, $\calo^{op}$ has $F(M)$ for its objects (for $M \in \calo$), and
induced homomorphisms for its morphisms. More generally, we can define $F
: \calh \to \calh$ in the same way.\\

Our analysis in the next few sections is in the spirit of \cite{BGG1},
\cite{GGOR}, and \cite{CPS1}.\\

\noindent {\bf Notation}: Throughout the rest of this paper (resp. in
the appendix), by the long exact sequence of Ext's, we mean the long
exact sequence of \Ext's in the abelian, self-dual category \caloo,
consisting of all objects of finite length in \calo (resp. in the abelian
category \calo). (That \caloo is abelian and self-dual will be proved
below.)

\begin{prop}\label{T2.1}
Let $M \in \calh$.
\begin{enumerate}
\item $ch_{F(M)} = ch_M$.
\item $F(F(M))$ is canonically isomorphic to $M$.
\item $\Hom(M,N) = \Hom(F(N),F(M))$ if $M,N \in \calh$.
\end{enumerate}
\end{prop}

The proof is standard, given that all weight spaces are
finite-dimensional, and hence reflexive.\\

\begin{prop}\label{Pduality}
\hfill
\begin{enumerate}
\item $F$ is an exact contravariant functor in \calh.
\item If $M \in \calh$ is simple, then so is $F(M)$. Further, $M = V(\la)
\Leftrightarrow F(M) = V(\la)$.
\item If $M \in \calo$ has a filtration in \calo with subquotients $V_i
\in \calo$, then $F(M)$ has a filtration in $\calo^{op}$, with
subquotients $F(V_i)$ occurring in {\em reverse order} to that of the
$V_i$'s.
\end{enumerate}
\end{prop}

\begin{proof}
We only show that if $M = V(\la)$, then $F(M) = V(\la)$. Now,
$\dim_k(F(M)_\la) = \dim_k(M_\la)$ (from Proposition \ref{T2.1}) = 1,
hence say $m^*$ spans $F(M)_\la$. Now, $m_\la \in M_\la$ is of maximal
weight, so $m^*$ is also maximal, and of weight $\la$. Therefore $Z(\la)
\to B_- m^* \to 0$, whence $0 \neq B_- m^* \subset F(M)$ simple. Thus,
$F(M) = V(\la)$.
\end{proof}

\begin{remark}
The last part is standard, once we verify that \calo is closed under
quotienting. Further, if $M \in \calo$ has finite length, then so does
$F(M)$, and $l(M) = l(F(M))$.\\
\end{remark}

$\calo$ is an additive category, with finite direct sums. All morphism
spaces are finite-dimensional. Inside \calo we define a new subcategory
\caloo, whose objects are all $M \in \calo$ with finite length (including
the zero module). Morphisms are module maps, as always.

\begin{theorem}\label{16finitelength}
\hfill
\begin{enumerate}
\item \calo is a full subcategory of $A$-mod, closed under taking
quotients.
\item In particular, every $M \in \caloo$ is a finite direct sum of
indecomposable objects.
\item \caloo is abelian, self-dual (i.e. $\caloo = \caloo^{op}$), and a
full subcategory of $A$-\emph{mod}.
\end{enumerate}
\end{theorem}

\begin{proof}[Sketch of proof]
For (1), if $M = \sum A m_i$, then $M/N = \sum A \overline{m_i}$, where
$0 \subset N \subset M$ is a submodule of $M \in \calo$. For (3), note
that if $M \in \caloo$ and $N$ is as above, then $l(N) \leq l(M) <
\infty$, so $N$ is finitely generated, and hence in \calo, thus in \caloo
as well. Thus \caloo is abelian. (This argument fails for \calo.)

To show that \caloo is self-dual, apply Proposition \ref{Pduality} above,
to any composition series for $M \in \caloo$.
\end{proof}\hfill

Inside \calo we have two sets of subcategories. For each $\la \in
\mf{h}^*$, we have the subcategory $\calo^{\leq \la}$ whose objects are
$M \in \calo$ so that $\Pi(M) \leq \la$. And for each $\bar{\la} \in
\mf{h}^* / (\Z \cdot \dd)$, we have the subcategory $\calo_{\bar{\la}}$,
whose objects are $M \in \calo$ so that $\Pi(M) \subset \la + \Z \cdot
\dd$.

\begin{prop}\label{split}
We work in the BGG category \calo.
\begin{enumerate}
\item $\calo^{\leq \la}$ is a full subcategory of $A$-mod, closed under
taking quotients.
\item If $N_\la = 0$ for some $N \in \calo$ and all $\la > \mu$, then
$\Ext^1(Z(\mu), N) = 0$.
\item If $Z(\mu) \to V \to 0$, then $\Ext^1(V,V(\mu))$ = 0.
\end{enumerate}
\end{prop}

\begin{proof}
(1) is easy to check, and the proof of (2) is as in \cite[Lemma
(16)]{Gu}. The proof of (3) is similar to that of (2), and we give it
below.\\

Say $0 \to V(\mu) \to M \mapdef{\pi} V \to 0$ is exact. Let $v_\mu$ be
the highest weight vector in $V$. Choose any (nonzero) $m \in
\pi^{-1}(v_\mu) \subset M_\mu$. Now, $v_\mu$ is maximal, so $\pi(N_+ m) =
0$, whence $N_+ m \subset V(\mu) \subset M$. But $V(\mu)$ has no weights
$> \mu$, so $N_+ m = 0$. Thus, $Z(\mu) \twoheadrightarrow B_- m
\overset{\pi}{\twoheadrightarrow} V$.

We know $m \notin V(\mu)$ because $\pi(m) \neq 0 = \pi(V(\mu))$. Now, say
$X = V(\mu) \cap B_- m$. Then $X$ is a submodule of $V(\mu)$ with
$\mu$-weight space zero, so $X=0$, and once more, we have $M = V(\mu)
\oplus B_- m$. So $B_- m \cong V$ and we are done.
\end{proof}

\begin{remark}
We cannot replace $Z(\mu)$ by a general $Z(\mu) \to V \to 0$ in part (2)
above, because we can have short exact sequences like $0 \to Z(\nu)
\hookrightarrow Z(\mu) \twoheadrightarrow Z(\mu) / Z(\nu) \to 0$. Also,
the above result says, in particular, that Verma modules and simple
modules have no self-extensions.\\
\end{remark}

\begin{prop}\label{extensions}
\hfill
\begin{enumerate}
\item If $Z(\la) \to N \to 0$ and $\Ext^1(Z(\mu),N) \neq 0$ (e.g. $N =
Z(\la), V(\la)$, etc.) then $\mu < \la$.
\item If $\Ext^1(V(\mu), V(\la)) \neq 0$ then it is finite-dimensional,
and $\mu < \la$ or $\la < \mu$.
\item Thus $\Ext^1(M,N)$ is finite-dimensional for $M,N \in \caloo$.
\end{enumerate}
\end{prop}

\begin{proof}
\hfill
\begin{enumerate}
\item This follows from the previous proposition: $\exists \omega > \mu$
so that $N_\omega \neq 0$. But since $N$ is standard cyclic, hence $\mu <
\omega \leq \la$, and we are done.\\

\item That $\mu \neq \la$ was shown in the previous proposition (since
there are no self-extensions). Now suppose $0 \to V(\la) \to M
\mapdef{\pi} V(\mu) \to 0$ is a nonsplit extension. The proof here is
similar in spirit to previous proofs. Say $v_\mu$ is the highest weight
vector in $V(\mu)$, and $m$ a lift to $M$. Then $\pi(N_+ m) = 0$, so we
have two cases.\\

\noindent $\bullet$ If $N_+ m = 0$ then $B_- m \twoheadrightarrow
V(\mu)$. Now, let $X = V(\la) \cap B_- m$, as earlier. $X$ is
nonzero since $M$ is a nontrivial extension, and so $X$ is a nonzero
submodule of $V(\la)$, whence $X = V(\la)$. But now $V(\la)
\hookrightarrow B_- m \twoheadrightarrow V(\mu)$, whence $\la < \mu$.\\

Now, since $X = V(\la)$, hence $\exists Z \in (B_-)_{\la - \mu}$ so
that $Zm = v_\la$ is the maximal vector in $V(\la)$. Conversely, any such
relation completely determines $M$, because $M_\mu$ is one-dimensional,
and $M$ has only two generators. Further, any such extension {\it has} to
be of this type, so $\dim_k(\Ext^1(V(\mu),V(\la))) \leq
\dim_k((B_-)_{\la - \mu}) = p(\la - \mu) < \infty$.\\

\noindent $\bullet$ If $N_+ m \neq 0$ then $(V(\la))_{\mu + \aaa} \neq 0$
for some $\aaa > 0$, whence $\mu < \mu + \aaa \leq \la$. But we are in
\caloo, because $M$ has length 2. Hence by Proposition \ref{extvecspc}
(in the appendix), $\Ext^1(V(\mu),V(\lambda)) \cong
\Ext^1(F(V(\lambda)),F(V(\mu))) = \Ext^1(V(\lambda),V(\mu))$ by
Proposition \ref{Pduality}, whence by the previous case it is
finite-dimensional.\\

\item  This follows from the previous part, using the long exact sequence
of \Ext's (and induction on lengths).
\end{enumerate}
\end{proof}

We now define the {\it co-standard modules} $A(\la) = F(Z(\la)) \in
\calo^{op}$. Since $Y(\la)$ was the radical of $Z(\la)$, and $V(\la)$ the
head, hence $V(\la)$ is the socle of $A(\la)$.\\

\section{Filtrations and finite length modules}

Note that to construct projectives in the classical case of \cite{BGG1},
one could quotient $\mf{Ug}$ by $(\mf{Ug}) \mf{n}_+^l$. Over here we
propose the following alternative:\\

Given $l \in \N$, look at the ``minimal weights" in $N_+^l$. That is,
define $\Sigma : (\Z \dd)^l \to \Z \dd$ by $(\mu_1, \dots \mu_l) \mapsto
\sum_{i=1}^l \mu_i$. Then the minimal weights in $N_+^l$ are simply $T =
\Sigma (\dd^l) = \{ \Sigma(i) : i \in \dd^l \}$. (Here, $\dd^l$ is the
$l$-fold Cartesian product of $\dd$.) Now define $B_{+l} = \sum_{\aaa \in
\nn \dd,\ \mu \in T} (B_+)_{\mu + \aaa}$.\\

\noindent Thus $\Pi(B_{+l})$ is closed under ``adding positive weights",
hence $B_{+l}$ is a two-sided ideal in $B_+$.\\

\noindent We claim that $B_+ / B_{+l}$ is finite-dimensional for all $l$.
Indeed, $\dd$ is finite, and any weight $\la$ of $B_+ / B_{+l}$ has to
look like $\sum_{\aaa \in \dd} c_\aaa \aaa$, where $0 \leq c_\aaa\
\forall \aaa$, and $\sum_\aaa c_\aaa < l$. Thus $\dim_k(B_+ / B_{+l})$ is
the sum of dimensions of finitely many weight spaces of $B_+$, each of
which is finite.\\

\noindent {\bf Definitions :}
\begin{enumerate}
\item Define the $A$-modules $P(\la,l)$ and $I(\la,l)$ (for $\la \in
\mf{h}^*$ and $l \in \N$) by
$$P(\la,l) = A / I_0(\la,l) \in \calo, \text{ and } I(\la,l) =
F(P(\la,l)) \in \calo^{op}$$
where $I_0(\la,l)$ is the left ideal generated by $B_{+l}$ and $\{(h -
\la(h) \cdot 1) : h \in \mf{h} \}$.
\item Given $\la \in \mf{h}^*$ and $l \in \N$, define the subcategory
$\calo(\la, l)$ to be the full subcategory of all $M \in \calo$ so that
$B_{+l} M_\la = 0$.
\item A (finite) filtration $0 = M_0 \subset M_1 \subset \dots \subset
M_n = M$ of an $A$-module $M$ is a\\
(a) {\it p-filtration} (cf \cite{BGG1}), denoted by $M \in \pfilt$, if
for each $i,\ M_i \in \calo$, and $M_{i+1} / M_i$ is a Verma module
$Z(\la_i)$.\\
(b) {\it q-filtration}, denoted $M \in \qfilt$, if for each $i,\ M_i \in
\calo^{op}$, and $M_{i+1} / M_i$ is a module of the form $A(\la_i)$.\\
(c) {\it SC-filtration} if for each $i,\ M_i \in \calo$, and $M_{i+1} /
M_i$ is standard cyclic.\\
\end{enumerate}

\noindent For example, if $l=1$, then we have $B_{+1} = N_+$, so
$P(\la,1) = Z(\la)$.

\begin{prop}\label{filtrn}
We still work in \calo.
\begin{enumerate}
\item Given $M \in \calo,\ M \in \pfilt$ iff $F(M) \in \qfilt$.
\item $\Hom(P(\la,l),M) = M_\la$ for each $M \in \calo(\la,l)$, so
$P(\la,l)$ is projective in $\calo(\la,l)$.
\item If $M \to N \to 0$ in \calo, and $M$ has an SC-filtration, then so
does $N$.
\item $P(\la,l)$ has an SC-filtration $\forall \la,l$.
\end{enumerate}
\end{prop}

\begin{proof}
\hfill
\begin{enumerate}
\item This follows from Proposition \ref{Pduality}, where we take each
$V_i$ to be a Verma module.\\

\item We know that $B_- \otimes H \otimes B_+ \twoheadrightarrow A
\twoheadrightarrow P(\la,l)$, and moreover, $H \otimes B_+
\twoheadrightarrow B_+ P(\la,l)_\la$. Therefore, $B_{+l}(H \otimes B_+)
\twoheadrightarrow B_{+l} P(\la,l)_\la$. Because $\mf{h}$ is
ad-semisimple, we see that $B_{+l}(H \otimes B_+) \subset A \cdot B_{+l}
\subset I_0(\la,l)$. Hence $B_{+l} P(\la,l)_\la = 0$, and $P(\la,l) \in
\calo(\la,l)$, as required.

Next, we show the exactness of $\Hom(P(\la,l),-)$. Given $\varphi \in
\Hom(P(\la,l),M)$, we get $v_\varphi = \varphi(1) \in M_\la$ (because $h
\varphi(1) = \varphi(h \cdot 1) = \la(h) \varphi(1)$ for each $h \in
\mf{h}$). Conversely, given $m \in M_\la$, define $\varphi \in
\hhh_k(k,M)$ by $\varphi(1) = m$. This extends to a map: $A \to M$ of
left $A$-modules. Because $M \in \calo(\la,l)$, hence $B_{+l}$ is in the
kernel, as is $(h - \la(h) \cdot 1)$. Thus $\varphi$ factors through a
map: $P(\la,l) \to M$ as desired. It is easy to see that both these
operations are inverses of each other, so we are done.\\

\item This is because quotients of standard cyclic modules are standard
cyclic.\\

\item The proof is similar to that in \cite{BGG1}. Moreover, the same
ordering holds among the terms of the filtration: if $Z(\la_{j+1}) \to
P_{j+1} / P_j \to 0$, and $\la_i \geq \la_j$, then $i \leq j$.
\end{enumerate}
\end{proof}

\begin{prop}\label{pfiltrn}
Suppose $M \in \pfilt$, and $S = \{ \nu \in \mf{h}^* :\ [M : Z(\nu)] \neq
0 \}$.
\begin{enumerate}
\item If $\la$ is maximal in $S$, then $\exists M'' \in \pfilt$ so that
$0 \to Z(\la) \to M \to M'' \to 0$ is exact.
\item If $\la$ is minimal in $S$, then $\exists M' \in \pfilt$ so that $0
\to M' \to M \to Z(\la) \to 0$ is exact.
\item Suppose $M_1, M_2 \in \caloo$. Then $M_1 \oplus M_2 \in \pfilt$ iff
$M_1,M_2 \in \pfilt$.
\end{enumerate}
\end{prop}

\begin{proof}
(1) and (3) follow from \cite{BGG1}, and (2) is cf.
\cite[(A3.1)(i)]{Don}.
\end{proof}\hfill

\noindent The next result comes from \cite{GGOR}, and involves
$\mf{h}$-diagonalizable modules $M$.

\begin{prop}\label{scfilt}
Suppose $M$ is $\mf{h}$-diagonalizable. Then the following are equivalent:
\begin{enumerate}
\item $M \in \calo$.
\item $M$ is a quotient of a direct sum of finitely many $P(\la,l)$'s.
\item $M$ has an SC-filtration. Further, the subquotients are standard
cyclic with highest weights $\la_i$, and we can arrange these so that
$\la_i \geq \la_j \Rightarrow i \leq j$.
\end{enumerate}
\end{prop}

\begin{proof}
We only show, in the part $(3) \Rightarrow (1)$, that $B_+$ acts locally
finitely on $M$. Since $M$ has an SC-filtration, $ch_M \leq \sum
ch_{Z(\la_i)}$, where we sum over a finite set. Thus, given $m \in
M_\mu$, we see that $\Pi(B_+ m) \subset \bigcup_i \{ \la : \mu \leq \la
\leq \la_i \}$, and each of these sets is finite. Hence $\Pi(B_+ m)$ is
finite, so $B_+ m$ is itself finite-dimensional.
\end{proof}\hfill

\begin{theorem}\label{moreoncalo}
Suppose every Verma module $Z(\la)$ has finite length.
\begin{enumerate}
\item Then \calo = \caloo.
\item If $\Ext^1(Z(\mu),M)$ or $\Ext^1(M,A(\mu))$ is nonzero for $M \in
\calo$, then $M$ has a composition factor $V(\la)$ with $\mu < \la$.
\item If $X \in \pfilt$ and $Y \in \qfilt$ then $\Ext^1(X,Y) = 0$.
\item If $X \in \pfilt$ and $Y \in \qfilt$ then
$$\dim_k(\Hom(X,Y)) = \sum_{\nu \in \mf{h}^*} [X : Z(\nu)][Y : A(\nu)]$$
where the terms on the RHS are the respective multiplicities in the
various filtrations. Thus
$$[X : Z(\mu)] = \dim_k(\Hom(X,A(\mu))), \emph{ and } [Y : A(\mu)] =
\dim_k(\Hom(Z(\mu),Y))$$
\end{enumerate}
\end{theorem}

\begin{proof}
\hfill
\begin{enumerate}
\item If all Verma modules have finite length, then so do all standard
cyclic modules, and since every module has an SC-filtration, hence all
modules have finite length.\\

\item This is cf. \cite[(A1.6)(ii)]{Don}.\\

\item The general case follows by the long exact sequence of \Ext's (and
induction on lengths of filtrations) from the case $X = Z(\mu),\ Y =
A(\la)$. To show the latter, suppose $\Ext^1(X,Y) \neq 0$. Applying the
previous part with $\Ext^1(Z(\mu),Y)$, we see that $Y$ has a composition
factor $V(\nu)$ with $\mu < \nu$. Since $Y = A(\la)$, hence we get $\mu <
\nu \leq \la$, so $\mu < \la$.

By symmetry, apply the previous part with $\Ext^1(X,A(\la))$, to get that
$X$ has a composition factor $V(\nu)$ with $\la < \nu$. Again, $\nu \leq
\mu$ because $X = Z(\mu)$, so $\la < \mu$. Thus we have obtained: $\la <
\mu < \la$, a contradiction. Hence all $\Ext^1(Z(\mu),A(\la)) = 0$.\\

\item For $X = Z(\mu), Y = A(\la)$, the result says that
$\dim_k(\Hom(Z(\mu),A(\la))) = \delta_{\mu \la}$, and this is simply
\cite[(A1.6)]{Don}. We again build the general case up, using the long
exact sequence of \Ext's and the previous part.
\end{enumerate}
\end{proof}

\section{Blocks in the BGG category $\calo$}\label{sec6}

Note that in the classical case, we had the notion of {\it blocks}
$\calo(\chi)$, where $\chi \in \hhh_{\C-alg} (\mf{Z}(\mf{U}(\mf{g})),
\C)$. Thus, a $\mf{g}$-module $V$ is in $\calo(\chi)$ iff for each $z$ in
the center $\mf{Z}$ one can find an $n$ so that $(z - \chi(z))^n$ kills
$V$. Furthermore, (cf. \cite[Exercise (23.9)]{H} or \cite[(7.4.8)]{Dix})
every algebra map from the center to $\C$ is of the form $\chi_\mu$ for
some $\mu \in \mf{h}^*$. Thus, the irreducible module $V = V(\la)$ is in
$\calo(\chi)$ iff $\chi_\la = \chi = \chi_\mu$, iff $\la + \delta$ and
$\mu + \delta$ are $W$-conjugate (by Harish-Chandra's theorem).

Over here, we do not have any of this, so we make some additional
assumptions. We make $\mf{h}^*$ into a (directed) graph as follows: given
$\la,\mu \in \mf{h}^*$, we say that $\la \to \mu$ if $Z(\la)$ has a
simple subquotient $V(\mu)$. Now make all edges non-directed, and for any
$\la \in \mf{h}^*$, define the set $S(\la) = \{ \mu : \la$ and $\mu$ are
in the same connected component of the graph $\mf{h}^* \}$.

\begin{stand}
$S(\la)$ is finite for each $\la$.\\
(Thus the $S(\la)$'s partition $\mf{h}^*$, and $S(\la) \subset \la + \Z
\dd$.)
\end{stand}

(For example, if $A = \mf{Ug}$, where $\mf{g}$ is a semisimple Lie
algebra over $\C$, then (it is well known that) the set $S(\la)$ is
contained in $W \bullet \la$, where the $\bullet$ denotes the twisted
action of the Weyl group: $w \bullet \la = w(\la + \delta) - \delta$,
where $\delta$ is the half-sum of positive roots.)\\

Note that category $\calo$ has the full subcategories $\calo(\la)$,
defined as follows: Given $\la \in \mf{h}^*,\ \calo(\la)$ contains
precisely those $M \in \caloo$, all of whose composition factors are of
the form $V(\mu)$, for some $\mu \in S(\la)$.

\begin{lemma}\label{caloo}
\calo = \caloo.
\end{lemma}

\begin{proof}
It suffices to show that every Verma module $Z(\la)$ has finite length.
Suppose $V$ is any subquotient of $Z(\la)$. Then $V$ has a maximal vector
$v_\mu$, so we get a nonzero module map $: A v_\mu = B_- v_\mu
\hookrightarrow V$. Hence $V(\mu) = A v_\mu / \Rad(A v_\mu)
\hookrightarrow V / \Rad(A v_\mu)$, so  $V(\mu)$ is a subquotient of
$Z(\la)$, and thus $\mu \in S(\la)$ by definition. We then claim that
$$l(Z(\la)) \leq \sum_{\mu \in S(\la)} \dim_k(Z(\la))_\mu = \sum_{\mu \in
S(\la)} p(\mu - \la) < \infty$$

because if $Z(\la) = V_0 \supset V_1 \supset \dots$, then each $V_i /
V_{i+1}$ has a maximal vector of weight $\mu$ for some $\mu \in S(\la)$.
Hence there can only be ``RHS-many" submodules in a chain, as claimed.
\end{proof}\hfill

\begin{theorem}\label{P5}
\hfill
\begin{enumerate}
\item $\Ext^1(V(\la'),V(\la)) = 0$ if $\la' \notin S(\la)$.
\item Given $M \in \calo$, let $S_M$ be the union of all $S(\la)$'s
corresponding to all simple subquotients of $M$. Suppose $S_M$ and
$S_{M'}$ are disjoint for $M,M' \in \calo$. Then $\hhom(M,M') =
\Ext^1(M,M') = 0$.
\item $\calo = \sum \calo(\la) = \bigoplus \calo(\la)$, where we sum over
all distinct blocks.
\end{enumerate}
\end{theorem}

\begin{proof}
\hfill
\begin{enumerate}
\item Say $0 \to V(\la) \to M \mapdef{\pi} V(\la') \to 0$ is a nontrivial
extension. Then we know from Proposition \ref{extensions} that $\la <
\la'$ or $\la' < \la$. Assume first that $\la' > \la$. Choose $m \in
\pi^{-1}(v_{\la'})$. Then from the proof of Proposition \ref{extensions},
we see that $V(\la) \hookrightarrow M = B_- m \twoheadrightarrow
V(\la')$. Hence $M$ is standard cyclic, so $Z(\la')$ has a simple
subquotient $V(\la)$, whence $\la' \in S(\la)$. On the other hand, if
$\la' \notin S(\la)$ and $\la > \la'$, then by Proposition
\ref{extvecspc} in the appendix, $\Ext^1(V(\la'),V(\la)) \cong
\Ext^1(V(\la),V(\la')) = 0$ (since $\calo = \caloo = \calo^{op}$).\\

\item This follows from (1) above, using induction on length, and the
long exact sequence of \Ext's. For the $\hhom$'s, use Corollary
\ref{C2.1} in place of part (1) above.\\ 

\item Given $M \in \calo$, we claim we can write it as $M = \bigoplus
M(\la)$, where $M(\la) \in \calo(\la)$. We prove this by using induction
on the length of $M$. For $l(M) =$ 0 or 1, we are easily done. Suppose we
have $0 \to N \to M \to V(\mu) \to 0$. We know that $N = \bigoplus
N(\la)$ because $N$ has lesser length.

Now $N = N' \oplus N(\mu)$, say, where $N'$ is the direct sum of all
other components of $N$. By Proposition \ref{extadditive} (in the
appendix), $M = N' \oplus M(\mu)$, where $0 \to N(\mu) \to M(\mu) \to
V(\mu) \to 0$. This is because $\Ext^1(V(\mu), N') = 0$ from the previous
part.\\

Thus $M = \bigoplus M(\la)$, where $M(\la) = M(\mu)$ if $\la = \mu$, and
$N(\la)$ otherwise.
\end{enumerate}
\end{proof}

\noindent{\bf Definition :}
Fix any indexing $S(\la) = \{ \la_1, \dots, \la_n \}$ that satisfies the
following condition: If $\la_i \geq \la_j$, then $i \leq j$. Now define
the \emph{decomposition matrix} $D$ in any block $\calo(\la)$ (where
$S(\la) = \{ \la_i \}$, under the above reordering) to be $D_{ij} =
[Z(\la_i) : V(\la_j)]$.

\begin{prop}\label{PDmatrix}
We work in a fixed block $\calo(\la)$.
\begin{enumerate}
\item $D$ is unipotent.
\item The Grothendieck group $Grot(\calo(\la))$ has the following
$\Z$-bases: $\{ [V(\mu)] :\ \mu \in S(\la) \},\ \{ [Z(\mu)] :\ \mu \in
S(\la) \},\ \{ [A(\mu)] :\ \mu \in S(\la) \}$.
\end{enumerate}
\end{prop}

\begin{remark}
Given $M \in \calo(\la)$, we now define the {\it multiplicities} $[M :
V(\la)],\ [M : Z(\la)] = [M : A(\la)]$ to be the coefficients of the
respective basis elements, when writing $[M]$ as a linear combination of
each of these bases. Then these actually equal the multiplicities of
$Z(\la)$'s and $V(\la)$'s in various p- and SC- filtrations (whenever $M$
does have such a filtration).\\
\end{remark}

\section{Projective modules in the blocks $\calo(\la)$}

Now fix $\la \in \mf{h}^*$. From above, we see that $\calo(\la)$ is a
full subcategory of \calo that is abelian, self-dual, and finite length.
We now construct projectives and progenerators in these blocks. Given
$\mu \in S(\la)$, as above we define $\calo(\la)^{\leq \mu}$ to be
$\calo(\la) \cap \calo^{\leq \mu}$.

\begin{prop}\label{projcover}
\hfill
\begin{enumerate}
\item If $V \in \calo(\la)^{\leq \la}$, then  $\Hom(Z(\la),V) \cong
V_\la$.
\item $Z(\mu)$ is the projective cover of $V(\mu)$ in $\calo(\la)^{\leq
\mu}$.
\end{enumerate}
\end{prop}

\begin{proof}
\hfill
\begin{enumerate}
\item We see that $\calo(\la)^{\leq \la} \subset \calo(\la,1)$, so
$P(\la,1) = Z(\la)$ is projective here.

\item We already know $Z(\mu)$ is an indecomposable projective in
$\calo(\la)^{\leq \mu} = \calo(\mu)^{\leq \mu}$, and $Y(\mu) =
\Rad(Z(\mu))$. Now use Theorem \ref{fitting} from the appendix.
\end{enumerate}
\end{proof}

\begin{theorem}\label{enoughproj}
\hfill
\begin{enumerate}
\item $\calo(\la)$ has enough projectives. 
\item There is a bijection between $S(\la)$ and each of the following
sets: indecomposable projectives (i.e. \emph{projective covers}),
indecomposable injectives (i.e. \emph{injective hulls}), Verma modules,
co-standard modules, and simple modules (all in $\calo(\la)$).
\item $\calo(\la)$ is equivalent to $(\emph{mod-}B_\la)^{fg}$, where
$B_\la$ is a finite-dimensional $k$-algebra.
\end{enumerate}
\end{theorem}

\begin{remark}
In fact, everything in Theorems \ref{fitting} and \ref{bass} holds here,
if we show the first part. For example, if $\la_0$ is maximal in
$S(\la)$, then $P(\la_0) = Z(\la_0)$ is the projective cover of
$V(\la_0)$, and $I(\la_0) = A(\la_0)$ is the injective hull.\\
\end{remark}

\begin{proof}
We only have to show that enough projectives exist in our abelian
category $\calo(\la)$. We refer to \cite[$\S 3.2 \S$]{BGS}. Following
Remark (3) there, we only need to verify five things (here) about
$\calo(\la)$, to conclude that enough projectives exist. We do so now.

\begin{enumerate}
\item $\mathcal{A} = \calo(\la)$ is a finite length abelian
$k$-category.\\

\item There are only finitely many simple isomorphism classes here
(because $S(\la)$ is finite).\\

\item Endomorphisms of any simple object (in fact, of any standard cyclic
object) are scalars, by Lemma \ref{Ltfae}.\\

The notation $\mathcal{A}_T$ refers precisely to $\calo(\la)^{\leq \mu}$. It
{\it is} a full subcategory. Further, $L(s) = V(s),\ \dd(s) = Z(s)$, and
$\nabla(s) = A(s)$ here. We also have maps $\dd(s) \to L(s)$ and $L(s)
\to \nabla(s)$.\\

\item As seen earlier, $Z(\mu) \to V(\mu)$ is a projective cover in
$\calo(\la)^{\leq \mu}$, and therefore $V(\mu) \to A(\mu)$ is an
injective hull, by duality. Both $Z(\mu)$ and $A(\mu)$ are
indecomposable, in particular.\\

\item $Y(s) = \ker(\dd(s) \to L(s))$ and $F(Y(s)) = \text{coker}(L(s) \to
\nabla(s))$ both lie in $\calo(\la)^{< s}$ for each $s \in S(\la)$
(meaning that they are in $\calo(\la)^{\leq s}$ and have no subquotients
$V(s)$).
\end{enumerate}
\end{proof}

\begin{remark}\hfill
\begin{enumerate}
\item The simple module, Verma module, co-standard module, projective
cover, and injective hull (of $V(\mu)$) corresponding to $\mu \in S(\la)$
are denoted respectively by $V(\mu),\ Z(\mu),\ A(\mu),\ P(\mu),\
I(\mu)$.

\item By duality, there are enough injectives in $\calo(\la)$. Since
$\calo = \bigoplus \calo(\la)$, hence \calo has enough projectives and
injectives; in particular, $P(\la)$ is projective and $I(\la)$ is
injective in \calo too. Every projective module $P \in \calo$ is of the
form $P = \bigoplus P(\la)^{\oplus n_\la}$, where only finitely many
$n_\la$'s are nonzero (and positive).\\
\end{enumerate}
\end{remark}

\noindent We conclude this section with one last result, cf. \cite{BGG1}.
It holds because $\calo = \bigoplus \calo(\la)$.

\begin{prop}\label{recip1}
Given $\la \in \mf{h}^*$ and $M \in \calo$, one has\\
$\dim_k(\Hom(P(\la),M)) = \dim_k(\Hom(M,I(\la))) = [M : V(\la)]$.\\
\end{prop}

\section{Every block $\calo(\la)$ is a highest weight
category}\label{sec8}

We now introduce the notion of a {\it highest weight category}, cf.
\cite{CPS1}, \cite[(A2.1)]{Don}. Let \calc be an abelian category over a
field $k$. Let $S$ index a complete collection of non-isomorphic simple
objects in \calc, say $\{ V(\la) : \la \in S \}$. We assume that \calc is
locally Artinian and satisfies the Grothendieck condition (these are
technical, though for our purposes, finite length would suffice), and
contains enough injectives.\\

The category \calc is then said to be a {\it highest weight category} if
$S$ satisfies the following conditions:

\begin{enumerate}
\item $S$ is an {\it interval finite} poset, i.e. there is a partial
ordering $\leq$ on $S$, and for each $\mu \leq \la \in S$, the set of
intermediate elements $[\mu,\la] = \{ \nu \in S : \mu \leq \nu \leq \la
\}$ is finite.\\

\item There is a collection of objects $\{ A(\la) : \la \in S \}$ of
\calc, and for each $\la$, an embedding $V(\la) \hookrightarrow A(\la)$,
such that all composition factors $V(\mu)$ of $A(\la) / V(\la)$ satisfy
$\mu < \la$. For $\mu,\la \in S$, we have that $\dim_k
\hhh_\calc(A(\la),A(\mu))$  and $\sup_{M \in J} [M : V(\mu)]$ are finite.
Here, $J$ is the set of all subobjects of $A(\la)$ of finite length, and
$[M : V(\mu)]$ denotes the multiplicity in $M$ of the simple module
$V(\mu)$.\\

\item Each simple $V(\la)$ has an injective envelope $I(\la)$ in \calc.
Further, the $I(\la)$'s each have a ``good filtration" which begins with
$A(\la)$ - namely, an increasing filtration $0 = F_0(\la) \subset
F_1(\la) \subset F_2(\la) \subset \dots$, such that:

\begin{enumerate}
\item $F_1(\la) \cong A(\la)$;
\item for $n>1$, $F_n(\la) / F_{n-1}(\la) \cong A(\mu)$ for some $\mu =
\mu(n) > \la$;
\item for a given $\mu \in S,\ \mu(n) = \mu$ for only finitely many $n$;
\item $\bigcup_i F_i(\la) = I(\la)$.\\
\end{enumerate}
\end{enumerate}

Reconciling this notation to our earlier notation, we see that each block
$\calc = \calo(\la)$ (is finite length, and hence) already satisfies all
conditions but two, namely, that $I(\la)/A(\la) \in \qfilt$, and each
co-standard cyclic factor $A(\mu)$ of $I(\la) / A(\la)$ satisfies $\mu >
\la$. (Here, we take $S$ to be the finite set $S(\la)$.)

\begin{stand}
The PBW theorem holds. In other words, $A \cong B_- \otimes_k H \otimes_k
B_+$.
\end{stand}

\noindent The final result in our analysis in this first part, is\\

\begin{theorem}\label{14projectives}
Every block $\calo(\la)$ is a highest weight category.
\end{theorem}

\noindent We need some intermediate results first.\\

\begin{prop}
\hfill
\begin{enumerate}
\item Fix $\la, \la' \in \mf{h}^*$. Then $\forall l \gg 0,\ \forall V \in
\calo(\la')$, we have $\Hom(P(\la,l),V) \cong V_\la$ as vector spaces.
\item $P(\la,l) \in \pfilt\ \forall \la,l$. Moreover, $[P(\la,l) :
Z(\la')] = p(\la - \la')$ if $\la' - \la \in \Pi(B_+ / B_{+l})$
(otherwise it is zero). Here $p$ is Kostant's function.
\item $P(\la) \in \pfilt$. If $[P(\la) : Z(\mu)] \neq 0$, then $\mu \geq
\la$.
\item $[P(\la) : Z(\la)] = 1$.
\end{enumerate}
\end{prop}

\begin{proof}
\hfill
\begin{enumerate}
\item The proof is similar to a proof in \cite{BGG1}.\\

\item Look at the analogous proof in \cite{BGG1}. Now that we know the
PBW theorem, that proof goes through completely.\\

\item Fix $l \gg 0$ so that $\Hom(P(\la,l),V) = V_\la$ for all $V \in
\calo(\la)$. Now suppose $P(\la,l) = \bigoplus_{\la'} N(\la')$. Since
$\Hom(P(\la,l),-)$ is exact in $\calo(\la)$, hence so is
$\Hom(N(\la),-)$. Thus $N(\la)$ is projective in $\calo(\la)$, so say
$N(\la) = \bigoplus_{\mu \in S(\la)} P(\mu)^{\oplus n_\mu}$.

Note that $\dim_k(\Hom(P(\la,l), V(\la))) = \dim_k(V(\la)_\la) = 1$, so\\
$\dim_k(\Hom(N(\la),V(\la))) = 1$ (because $\calo = \bigoplus
\calo(\la)$). Applying Proposition \ref{recip1}, we get $n_\la = 1$. Thus
$P(\la)$ is a direct summand of $P(\la,l)$, and $P(\la,l)$ has a
p-filtration, so by Proposition \ref{pfiltrn}, $P(\la) \in \pfilt$.

Finally, $P(\la)$ is a summand of $P(\la,l)$, hence {\it for all} $\mu$
we have $[P(\la) : Z(\mu)] \leq [P(\la,l) : Z(\mu)] \leq p(\la - \mu)$.
Therefore $[P(\la) : Z(\mu)] \neq 0$ only if $\la \leq \mu$.\\

\item Suppose $P(\la) \supset M_1 \supset \dots$ is a p-filtration, with
$P(\la) / M_1 \cong Z(\mu)$ for some $\mu \geq \la$. Then $P(\la)
\twoheadrightarrow P(\la) / M_1 = Z(\mu) \twoheadrightarrow Z(\mu) /
Y(\mu) = V(\mu)$ simple. Hence the composite has kernel $\Rad(P(\la))$,
whence $V(\mu)= V(\la)$, or $\mu = \la$. Hence $[P(\la) : Z(\la)] > 0$.
Also, $[P(\la) : Z(\la)] \leq [P(\la,l) : Z(\la)] = p(\la - \la) = 1$, so
we are done.
\end{enumerate}
\end{proof}

\begin{proof}[Proof of the Theorem]
Dualize the p-filtration for $P(\la)$ (in the last part above) to get a
q-filtration for $I(\la)$. Clearly, $P(\la) / M_1 = Z(\la)$ means that
the filtration looks like $0 \subset A(\la) \subset \dots \subset
I(\la)$. The weights are suitably ordered, hence $\calo(\la)$ is a
highest weight category.
\end{proof}\hfill

From above, we conclude that every projective module in $\calo$ has a
p-filtration, since each $P(\la)$ does. Also, since $\calo(\la)$ is a
highest weight category, we have {\it Brauer-Humphreys / BGG
Reciprocity}, which says that $[P(\la) : Z(\mu)] = [A(\mu) : V(\la)] =
[Z(\mu) : V(\la)]$. Further, the cohomological dimension of $\calo(\la)$
is bounded above, hence finite.

There are many more results, especially on Tilting modules and Ringel
duality, which are readily found in \cite{Don}, for instance, and which
we do not mention here.\\

\section*{\bf Part 2 : The (deformed) symplectic oscillator algebra
$H_f$}

In this part, we show that all assumptions in the first part are true for
the algebra $H_f$, which we shall define presently. We prove the PBW
theorem for $H_f$, classify all finite-dimensional simple modules, state
the well-known character formulae, and take a closer look at Verma
modules. We conclude by producing a counterexample to Weyl's theorem (of
complete reducibility) for a special case.\\

\section{Introduction; automorphisms and anti-involutions}

\noindent We continue to work over an arbitrary field $k$ of
characteristic zero.\\

Consider the Lie algebra \spn. The Cartan subalgebra $\mf{h}$ has basis
$h_i = e_{ii} - e_{i+n, i+n} \ (1 \leq i \leq n)$, though these do not
correspond to the simple roots of \spn. Now define the functionals
$\eta_i \in \mf{h}^*$ by $\eta_i(h_j) = \delta_{ij}$. Then the roots and
root vectors are:

\noindent $u_{jk} = e_{jk} - e_{k+n,j+n} : 1 \leq j \neq k \leq n$ (root
= $\eta_j - \eta_k$)\\
$v_{jk} = e_{j,k+n} + e_{k,j+n} : 1 \leq j < k \leq n$ (root = $\eta_j +
\eta_k$)\\
$w_{jk} = e_{j+n,k} + e_{k+n,j} : 1 \leq j < k \leq n$ (root = $-\eta_j -
\eta_k$)\\
$e_j = e_{j,j+n} : 1 \leq j \leq n$ (root = $2 \eta_j$)\\
$f_j = e_{j+n,j} : 1 \leq j \leq n$ (root = $-2 \eta_j$)\\

\noindent The simple roots are given by $\{ \eta_i - \eta_{i+1} : 0 < i <
n \}$ and $2 \eta_n$.

\begin{remark}\label{R7}
It is easier for calculations to use $e_j = 2 e_{j,j+n}$ and $f_j =
2e_{j+n,j}$, because then $h_j = u_{jj}, e_j = v_{jj}, f_j = w_{jj}$.\\
\end{remark}

\noindent Let $B = k[X_1, \dots, X_n]$, and consider a $2n$-dimensional
$k$-vector space $V \subset \text{End}(B)$, with basis given by $\{ X_i =
\text{multiplication by } X_i : 1 \leq i \leq n,\ Y_i = (\partial /
\partial X_i) : 1 \leq i \leq n\}$. Then the subalgebra generated by $V$
in End($B$) is called the {\it Weyl algebra} = $A_n$. We now construct
the {\it Weil representation} of \spn on $B$. More precisely, define the
map $\varphi : \uu \to A_n \subset \mf{gl}(B)$ as follows:
$$h_i \mapsto X_i Y_i + 1/2, \quad u_{jk} \mapsto X_j Y_k, \quad v_{jk}
\mapsto - X_j X_k, \quad w_{jk} \mapsto Y_j Y_k, \quad e_j \mapsto -X_j^2
/ 2, \quad f_j \mapsto Y_j^2 / 2$$

Thus we obtain a representation $\varphi_0 : H_0 \to A_n$, where $H_0 =
\uu \ltimes A_n$, and $\varphi_0 = \varphi \ltimes id$. (It is a {\it
faithful} map of Lie algebras: $\spn \to A_n$.) Here $H_0$ is defined by
$Za - aZ = Z(a)\ (= [\varphi(Z),a])$, where $Z \in \spn,\ a \in V$, and
$Z(a)$ is the action of $Z$ on $a$. Thanks to our choice of $\varphi$,
this also agrees with the natural action of $\mf{sp}(2n)$ on $V$ (i.e. as
$2n \times 2n$ matrices, acting on vectors in $V$).

Note that $H_0$ arises from the {\it symplectic oscillator algebra}
$\mf{sp}(2n) \ltimes \mf{h}_n$ (relations as above) by: $H_0 =
\mf{U}(\mf{sp}(2n) \ltimes \mf{h}_n) / (I - 1)$, where $I$ is the central
element in (the $(2n+1)$-dimensional {\it Heisenberg} algebra)
$\mf{h}_n$.\\

We now consider a deformation over $k[T]$ of $H_0$. For $f \in k[T]$,
define $H_f = T(V_0) / \langle R_f \rangle$, where $V_0 = \spn \oplus V$
and $R_f$ is generated by $Za - aZ = Z(a)$, the usual $\mf{sp}(2n)$
relations, $[X_i, X_j],\ [Y_i, Y_j]$, and the deformed relations $[Y_i,
X_j] - \delta_{ij}(1+f(\dd))$. Here, $\dd$ is the quadratic Casimir
element in \spn, acting on $A_n$ via the above map $\varphi$, as the
scalar $c_\varphi = -(2n^2+n) / 16(n+1) \in \mathbb{Q} \subset k$.

\begin{remark}
For $n=1$, we can show that \spn commutes with all of $[V,V]$, so that
the deformation must lie in $\mf{Z}(\mf{U}(\mf{sp}(2n)))$, and for $n=1$,
this is precisely $\C[\dd]$. This explains the choice of deformed
relations. (However, $\dd$ does not commute with all of $V$.)
\end{remark}

\noindent We now explicitly describe some automorphisms and an
anti-involution of $H_f$.\\

\noindent {\it Anti-involutions :}
Define $i : V_0 \to V_0$ by sending $X_j \mapsto Y_j, Y_j \mapsto X_j,
u_{jk} \mapsto u_{kj}, v_{jk} \mapsto -w_{jk}, w_{jk} \mapsto -v_{jk}\
\forall j,k$ (as in Remark \ref{R7} above). This extends to an
anti-involution : $T(V_0) \twoheadrightarrow H_f$, defined on monomials
by reversing the order, and this map {\it does} vanish on $R_f$, as
desired. In addition, it takes $\mf{U}(N_+)_\mu$ to $\mf{U}(N_-)_{-\mu}$
for every $\mu$, and acts on $\mf{h}$ as the identity.\\

\noindent {\it Automorphisms / lifts of the Weyl Group :} Let us now lift
the Weyl group to automorphisms of $H_f$. Let $S = \{ u_{jk}, v_{jk},
w_{jk}, X_j, Y_j \}$. Then $\forall a_\aaa \in S \cap (\spn)_\aaa$, we
see that $\tau_{a_\aaa}(b) := \exp(\text{ad } a_\aaa)(b)$ is a finite
series $\forall b \in S$, if $\aaa \neq 0$. Further, $\tau_\aaa :=
\tau_{a_\aaa} \tau_{-a_{-\aaa}} \tau_{a_\aaa}$ takes $(V_0)_\mu$ to
$(V_0)_{\sigma_\aaa(\mu)}$ for all (simple) roots $\aaa$. In addition, it
also permutes the Cartan subalgebra $\mf{h}$ ``appropriately". Thus each
$\tau_\aaa$ is an algebra automorphism, preserving $V_0$ and taking
$(H_f)_\mu$ to $(H_f)_\nu$, where $\nu = \sigma_\aaa(\mu)$.

Now, we know (cf. \cite[Exercise (13.5)]{H}), that the Weyl group $W =
(\Z / 2 \Z)^n \rtimes S_n$ of \spn contains $-1$. So we can construct an
automorphism $\tau$ of $H_f$ that restricts to $-1$ on $\mf{h}$,
preserves $V_0$, and takes each weight space to the corresponding
negative weight space.\\

\section{Standard cyclic $H_f$-modules in the BGG category}

Let $\Phi$ (resp. $\Phi_f$) be the root system of $\spn$ (resp. $H_f$).
Then $\Phi_f = \Phi \coprod \{ \eta_i, - \eta_i : 1 \leq i \leq n\}$, and
$\dn = 1+f(\dd)$. We write positive and negative roots as $\Phi_f^+ =
\Phi^+ \coprod \{ \eta_j \}$ and $\Phi_f^- = -\Phi_f^+$. Similar to
\cite{H}, we introduce an ordering among the roots as follows: $\la
\succ_f \mu$ if $\la - \mu$ is of the form $(m \eta_n + \sum_{i<n} k_i
\alpha_i)$, where $m,k_i \in \nn$, and $\alpha_i = \eta_i - \eta_{i+1}$
are the first $n-1$ simple roots (as above).

Now define Lie subalgebras $N_+ = [B_+, B_+] \subset B_+ \subset H_f$ as
follows: $B_+ = \mf{h} \bigoplus N_+$ is a Borel subalgebra, and $N_+ =
\bigoplus_{i=1}^n kX_i \oplus \bigoplus_{\alpha \in \Phi^+}
(\spn)_\alpha$ is nilpotent. Similarly, we have $B_-$ and $N_-$. (Note
that these are {\it not} the $B_\pm,N_\pm$ of Section \ref{sec5} above;
rather, those are given here by $\mf{U}(B_\pm),\mf{U}(N_\pm)$.)\\

We now observe that the ``Setup" for the analysis in the first part of
this paper is partially valid here. The assumptions in Section \ref{sec5}
are all satisfied. Thus Theorem \ref{2tfae} holds here. Assuming the PBW
theorem, we introduce another equivalent condition:

\begin{cor}\label{C5}
Suppose $ H_f \cong \mf{U}(N_-) \otimes_k \mf{U}(\mf{h}) \otimes_k
\mf{U}(N_+)$. Then all nonzero maps from $Z(\mu)$ to $Z(\la)$ are
injections.
\end{cor}

The proof uses the fact that $\mf{Ug}$ is an integral domain for {\it
any} Lie algebra $\mf{g}$ (cf. \cite[(2.3.9)]{Dix}).\\

Now suppose $V(\la)$ is finite-dimensional. Since any $H_f$-module is
also a \spn-module, hence Weyl's theorem applies (cf. \cite[$\S 7.8
\S$]{We}), and $V(\la)$ is a direct sum of finitely many $V_C(\mu)$'s,
where $V_C(\mu)$ is the irreducible \spn-module of highest weight $\mu$
(which is dominant integral because $V$ has finite dimension). Thus if
$V(\la)$ is finite-dimensional, then $\la \in \Lambda^+$. Further,
$\Pi(V)$ is saturated (under the action of the Weyl group $W$ of \spn).

We now come to character theory. $W$ acts naturally on $\Z[\La]$ by
$\sigma e(\la) = e(\sigma \la)$. If $\dim_k(V) < \infty$, then
$\dim(V_\mu) = \dim(V_{\sigma (\mu)})$, i.e. $ch_V \in \Z[\La]^W$. Let us
define $\tau_\aaa \in \text{Aut}(V)$ for any finite-dimensional module
$V$. Since all nonzero root vectors in \spn act nilpotently on $V$, we
can define $\tau_\aaa$ as above. Then $\tau_\aaa \in \text{Aut} (V)$ and
$\tau_\aaa : V_\mu \to V_{\sigma_\aaa(\mu)}$ by $\mf{sp}(2n)$-theory. In
particular, we again get $ch_V \in \Z[\La]^W$.

In order to handle infinite-dimensional modules, we redefine the formal
character as a function $: \La \to \Z$. Then multiplication becomes
convolution. The $e(\mu)$ becomes $\epsilon_\mu : \nu \mapsto \delta_{\mu
\nu}$, so $\sigma (\epsilon_\mu) = \epsilon_{\sigma \mu}$. The usual
definition of the Kostant function now coincides with our previous
definition (setting $B_- = \mf{U}(N_-)$). The {\it Weyl function} $q$ is
just $\prod_{\aaa \in \Phi_f^+} (e(\aaa /2) - e( -\aaa /2))$, and we set
$\delta = \frac{1}{2} \sum_{\aaa \in \Phi_f^+} \aaa$.\\

\begin{lemma}\label{L3.3} Assume the PBW theorem holds. Then
\begin{enumerate}
\item $p = ch_{Z(0)}$
\item $ch_{Z(\la)} = p * \epsilon_\la$
\item $q * ch_{Z(\la)} = q * (p * \epsilon_\la) = \epsilon_{\la +
\delta}$.
\end{enumerate}
\end{lemma}

\noindent The proof is a matter of easy calculation.\\

\section{$H_f$-modules for $n=1$}\label{hfn1}

Throughout the rest of this paper, we take $n=1$. Thus our Lie algebra is
$C_1 = \spl = \mf{sp}(2)$. We denote the generators of $H_f$ by
$E,F,H,X,Y$. The ``root system" is $\Phi_f = \{ \pm \eta,\ \pm 2 \eta
\}$, and the Weyl group $W$ is simply $S_2$. We may also prefer to work
with a related group $W' = S_2 \times S_2$, whose action on the weights
will be seen later, in $\S \ref{sec19} \S$ below.

We write down the generators and relations explicitly here. $H_f$ is
generated by $X,Y,E,F,H$, with $E,F,H$ spanning \spl. The other relations
are: $[E,X] = [F,Y] = 0,\ [E,Y] = X,\ [F,X] = Y$. Further, $X$ and $Y$
are weight vectors for $H$ : $[H,X] = X,\ [H,Y] = -Y$. Finally, the
deformed relation is $[Y,X] = \dn = 1 + f(\dd)$, where $\dd$ is the
quadratic Casimir element $\frac{1}{4}(EF+FE+H^2/2)$.

Note that the original symplectic oscillator algebra contains the {\it
oscillator algebra} \cala (cf. \cite{KalMil}), where $E_+ = X,\ E_- = Y,\
H = H,\ \mathcal{E} = I = 1$ (where $I$ is the central element in
$\mf{h}_1$).\\

Our main motivation is to prove the PBW theorem, and the remaining
``standing assumption" mentioned in Section \ref{sec6} above (note that
all Verma modules are automatically nonzero if PBW holds). However, we
will also consider other things - for example, the structure of
finite-dimensional modules and Verma modules.

First of all, notice (cf. \cite{H}) that on any standard cyclic
\usl-module, $\dd$ acts by a scalar. Therefore $\dn$ also acts by a
scalar, and let us denote this by $c_{0r}$ if the module is of highest
weight $r \in k$. Clearly, $c_{0r}$ depends on the polynomial $f$ as
well.\\

We now come to calculations. First of all, observe that $\mf{U}(N_-) =
k[Y,F]$ because $YF = FY$. Thus we see that in $Z(r)$, a spanning set for
the $(r-m)$-weight space is $Y^m,\ Y^{m-2}F,\ \dots$. Define the
constants
\begin{equation}\label{E1}
\aaa_{rm} = \sum_{i=0}^{m-2} (r+1-i)c_{0,r-i} \quad \mbox{and } d_{r-m}
= \aaa_{rm} / (r-m+2)(r-m+3)
\end{equation}

Of course, to define $d_{r-m}$ we should not have $r = m+2,\ m+3$. Also,
we clearly have $m \in \N$ (for $m=1$ we can take the empty sum = 0).\\

For the time being, we work only with standard cyclic modules. Consider
any $Z(r) \to V = H_f v_r \to 0$, for $r \in k$. We have

\begin{theorem}\label{3basics}
Let $V = H_f v_r$. Then
\begin{enumerate}
\item $v_r$ and $v_{r-1} = Yv_r$ are \spl-maximal vectors (i.e. $Ev_r =
Ev_{r-1} = 0$).\\

Now say $t \in r - 2 - \nn$. Wherever $d_t$ can be defined, we have $R_t$
and define $S_t$ inductively:

\begin{equation}\label{Rt}
Xv_{t+1} = EY v_{t+1} = -\frac{\aaa_{r, r-t}}{t+3} v_{t+2}
\tag{$R_t$}
\end{equation}

\begin{equation}\label{St}
v_t \overset{\mbox{{\em def}}}{\mbox{ = }} Yv_{t+1} + d_t Fv_{t+2}
\tag{$S_t$}
\end{equation}

\noindent For the same values of $t$, we also have the following:\\

\item $v_t = p_{r-t}(Y,F)v_r$ for some polynomial $p_{r-t}(Y,F) = Y^{r-t}
+ c_1 FY^{r-t-2} + \dots \in k[Y,F]$ (monic in $Y$).
\item Say $v \in V_t$. Then $Ev = 0$ iff $v \in k \cdot v_t$.\\
\end{enumerate}
\end{theorem}

\begin{remark}\hfill
\begin{enumerate}
\item Thus, if $r \in \nn$, then the equations are valid until we reach
$t=-1$. We can define $v_{-1}$ and calculate $Xv_{-1}$, but cannot go
beyond that. Of course, if $r \notin \nn$ then we can go on
indefinitely.\\

\item Suppose $t > -2$ or $t \notin \nn$. Then we can rewrite $(R_t)$ as
\begin{equation}\label{Rt'}
Xv_{t+1} = EY v_{t+1} = -(t+2)d_t v_{t+2}
\tag{$R_t'$}
\end{equation}

\item Henceforth, the phrase ``where(ver) $d_t$ can be defined" means
``where(ver) $t>-2$ if $r \in \nn$".\\
\end{enumerate}
\end{remark}

\begin{proof}[Proof of the theorem]
This is just inductive calculations.
\end{proof}

\begin{cor}\label{C6}
Suppose $v_t, v_{t+1} \neq 0$ for some $t\ (t > -2$ if $r \in \nn)$. Then
$v_t$ is maximal iff $\aaa_{r,r-t+1} = 0$.
\end{cor}

We will see further below that one implication holds for any $r \in k$,
namely, that if $v_t$ is maximal in $Z(r)$, then $\aaa_{r,r-t+1} = 0$.

\begin{cor}\label{C7}
Suppose $v_t = 0$. If $v_{t-n}$ can be defined for $n \in \nn$, then
$v_{t-n}=0$.\\
\end{cor}

\begin{cor}\label{C9}
Suppose $V$ (as above) has another maximal vector $v_t$ for some $t \in r
- \N$. Then a weight vector $v_T$ in $V' = H_f v_t$ (defined in $V'$ by
the relation $(S_T)$ for some $T$, so that $d_{T-1}$ is defined) is
maximal in $V'$ iff it is maximal in $V$.
\end{cor}

\begin{proof}
The proof is, of course, that a maximal vector generates a submodule, and
a submodule of a submodule is still a submodule. However, there is a
related phenomenon occurring among the $\aaa_{rm}$'s. The point is that
if $H_f v_T \subset H_f v_t \subset H_f v_r = V$ are all submodules of
$V$, then these $v$'s are maximal vectors, and Corollary \ref{C6} says
that there is a relation among the various $\aaa_{rm}$'s. In fact, it is
easy to show (from definitions) that
\begin{equation}\label{eqn1}
\aaa_{r,r-T+1} = \aaa_{r,r-t+1} + \aaa_{t,t-T+1}
\end{equation}
\end{proof}

\begin{cor}\label{C10}
Say $V = V(r)$ is simple, and $d_t$ can be defined for $t \in r - 2 -
\nn$. Then $d_{t-1}=0$ only if $v_t=0$.\\
\end{cor}

\section{General philosophy behind the structure theory}

As we shall see, many standard cyclic (resp. Verma, simple) $H_f$-modules
$Z(r) \to V \to 0$, are a direct sum of a progression of standard cyclic
(resp. Verma, simple) \usl-modules $V_{C,t}$ of highest weight $t \in r
- \nn$. (Each module $V_{C,t}$ has multiplicity one as well.)

If this progression terminates, say at $Z_C(t) \to V_{C,t} \to 0$ for
some $t = r-n$, then (we show later that) $\alpha_{r,n+1} = 0$. The
converse is true, for instance, when $r \notin \nn$ (as the results and
remarks in the previous section suggest), or if $V$ is finite-dimensional
simple (as we shall see in a later section). But there are
counterexamples to a general claim of this kind, which we shall provide
below.

The specific equations governing such a direct sum $V = \oplus_i
V_{C,r-i}$ are the subject of the previous subsection. Very briefly,
though, if $v_t$ is the highest weight vector (for \usl) in $V_{C,t}$,
then we see that $E(Xv_t) = X(Ev_t) = 0$, so that $Xv_t$ must be a
highest weight vector in $V_{C,t+1}$. Since the highest weight space in
each $V_{C,t}$ is one-dimensional, there is some scalar $a_t$ so that
$Xv_t = a_t v_{t+1}$. And if {\it this} scalar vanishes, then $v_t$ is
$H_f$-maximal in $V$.\\

{\it This} is the scalar $\alpha_{r,n}$ (upto a constant).\\

\section{Certain Verma modules are nonzero}\label{S11}

We now show that $Z(r)$ is nonzero if $r \notin \nn$. In fact, we show it
to be isomorphic to $\mf{U}(N_-)$, by constructing a standard cyclic
module of highest weight $r$, whose character is $ch_{\mf{U}(N_-)} *
\epsilon_r$.

\begin{lemma}\label{Lbasis1}
We work in $H_f$.
\begin{enumerate}
\item $\displaystyle [X, F^jY^i] = - F^j \sum_{l=0}^{i-1} Y^{i-l-1} \dn
Y^l - j F^{j-1} Y^{i+1}$
\vspace{1ex}
\item $[E, F^jY^i] = - \displaystyle F^j \sum_{m=0}^{i-2} (i-1-m)
Y^{i-2-m} \dn Y^m + j(r-i-j+1) F^{j-1}Y^i$
\end{enumerate}
\end{lemma}

\begin{proof}
We show by induction that $[F^j,X] = jF^{j-1}Y$. Now the proof is just
small calculations.
\end{proof}\hfill

Now fix $r \notin \nn$. Define a module $V$ with $k$-basis $\{ v_{ij} :
i,j \in \nn \}$. We now define the module structure by: $Yv_{ij} =
v_{i+1,j},\ Fv_{ij} = v_{i,j+1},\ Hv_{ij} = (r-i-2j)v_{ij}$. For the $E$-
and $X$-actions, we use the preceding lemma as follows:\\

We first set $Xv_{00} = Ev_{00} = 0$. From above, $Y^k F^l v_{ij} =
v_{i+k,j+l}$, so $YF = FY$ (on all of $V$). Now we multiply both sides of
the equations in the lemma above, by $v_{00}$ on the {\it right}. The
left hand sides give us $Xv_{ij}$ and $Ev_{ij}$ respectively. The right
hand sides are calculated inductively, starting from the fact that we set
$Xv_{00} = Ev_{00} = 0$. We see that we can define $\dd v_{ij}$
inductively, using the above lemma; hence we can also define $\dn v_{ij}$
using induction on $(i,j)$.\\

This is how we define $Xv_{ij}$ and $Ev_{ij}$ inductively. Now we need to
verify that the module structure is consistent with the relations in
$H_f$. (To start with, it is easy to compute that $Ev_{10} = EYv_{00} =
0$. Similarly, $\dn v_{00} = c_{0r} v_{00}$ and $\dn v_{10} = c_{0,r-1}
v_{10}$.)

First of all, one sees from above that the $E,X,Y,F,H$-actions take
weight vectors into appropriate weight spaces, so all relations of the
form $[H,a_\mu] = \mu(H) a_\mu$ automatically hold. As seen above, $YF =
FY$. We now verify the following:
$$[E,Y] = X \qquad [F,X] = Y \qquad [E,F] = H \qquad [Y,X] = \dn$$

Let us show that $EY-YE = X$; the others are similar (and easy). Note
that in the calculations below, the right hand side quantities are to be
(right) multiplied by $v_{00}$.
\begin{eqnarray*}
EY v_{ij} & = & -F^j \sum_{m=0}^{i-1} (i-m) Y^{i-1-m} \dn Y^m + j(r-i-j)
	F^{j-1} Y^{i+1}\\
YE v_{ij} & = & -F^j \sum_{m=0}^{i-2} (i-1-m) Y^{i-1-m} \dn Y^m +
	j(r-i-j+1) F^{j-1} Y^{i+1}\\
X v_{ij} & = & -F^j \sum_{m=0}^{i-1} Y^{i-1-m} \dn Y^m -j F^{j-1}Y^{i+1}
\end{eqnarray*}

\noindent To verify the last relation, namely $EX=XE$, we now introduce
another basis of $V$.

\begin{lemma}\label{Lbasis2}
The set $\{ F^j v_{r-n} : j,n \in \nn \}$ is a basis for $V$, where
$v_{r-n}$ is defined in equation \eqref{St}.
\end{lemma}

\begin{proof}
The equations \eqref{Rt},\eqref{St} hold for all $t = r-n$ (since $r
\notin \nn$), so define (for all $n$) $v_t = v_{r-n} = p_n(Y,F) v_r$,
where all $p_n$'s are {\it monic}. This makes a change of basis easy to
carry out.
\end{proof}\hfill

\begin{remark}
Until now, we have {\it never} used the relation $EX=XE$. We now define
some module relations using the $F^j v_{r-n}$'s. That they hold can be
checked from the relations \eqref{Rt} and \eqref{St}, once again without
using $[E,X]=0$.
\end{remark}

\noindent $H \cdot F^j v_{r-n} = (r-n-2j)F^j v_{r-n}\\
E \cdot F^j v_{r-n} = j(r-n-j+1) F^{j-1}v_{r-n} \\
X \cdot F^j v_{r-n} = -jYF^{j-1} v_{r-n} - (r-n+1) d_{r-n-1} F^j
v_{r-n+1}$ (Here, $d_{r-1} = 0$ as above.)\\

We now verify the remaining relation, namely, $EX=XE$. Note that we are
free to use the other relations now, since we showed above that they hold
on all of $V$. We compute\\
$EX (F^j v_{r-n}) = -j(r-n-j+1)[ (j-1)YF^{j-2}v_{r-n} + (r-n+1) d_{r-n-1}
F^{j-1} v_{r-n+1} ] = XE (F^j v_{r-1})$.\\

We have thus checked all relations, and hence shown that there exists a
nonzero standard cyclic module $Z(r) \to V \to 0$ of highest weight $r
\notin \nn$. In fact,

\begin{theorem}\label{Tverma}
$0 \neq Z(r) \cong k[Y,F]\ \forall r \notin \nn$.\\
\end{theorem}

\section{$\aaa_{rm}$ is a polynomial}

We now show that $\aaa_{rm}$ is a polynomial in two variables. Actually
we show a more general result, that can be applied to various
``polynomials" in our setting. Throughout, by deg($f$) we mean the degree
of $1+f(T)$, because that is what we use in handling $\dn$.

\begin{prop}\label{7sumofpowers}
Given $d \in \nn$, there exists a polynomial $g_d \in \mathbb{Q}[T]
\subset k[T]$, of degree $d+1$, so that $g_d(0) = 0$, and $g_d(T) -
g_d(T-1) = T^d$.
\end{prop}

\begin{proof}
We inductively define $\displaystyle g_d(T) = \frac{1}{d+1} \bigg[
(T+1)^{d+1} - 1 - \sum_{i=0}^{d-1} \binom{d+1}{i} g_i(T) \bigg]$. The
base case is $g_0(T) = T$. Then one checks that $g_d$ is as desired, by
induction on $d$. (In particular, for all $m \in \nn$, we have $g_d(m) =
\sum_{n=1}^m n^d$, e.g. $g_1(T) = T(T+1)/2$.)
\end{proof}\hfill

\begin{cor}\label{C11}
$\aaa_{rm}$ is a polynomial in $r,m$, of degree $2 \deg(f) + 2$ in $m$,
and degree $2\deg(f) + 1$ in $r$.
\end{cor}

\begin{proof}
First of all we find out what $c_{0r}$ actually is - or more precisely,
what $\dd$ acts on $\usl v_t$ by. So suppose we have $E v_t = 0$. Then
$\dd = (EF + FE + H^2 / 2) / 4$ acts on $v_t$ by: $(EF v_t + F.0 + H^2
v_t / 2) / 4 = (tv_t + 0 + t^2 v_t / 2)/4 = [(t^2 + 2t)/8] v_t$.

Thus $\dd$ acts on $v_t$ by the scalar $c_t = (t^2+2t)/8$. Remember, of
course, that $t$ is of the form $r-m$ for some $m \in \nn$. Now, we see
that $\dn$ acts on $v_t$ by $c_{0,r-m} = 1+f(c_{r-m})$. This is clearly a
polynomial in $r$ and $m$, if we expand out $f(c_{r-m})$ formally.

Now equation \eqref{E1}, combined with Proposition \ref{7sumofpowers},
says that $\aaa_{r,m}$ is a polynomial in two variables, as required.
Also, $1+ f(c_t)$ is of degree 2 deg($f$) in each of $r$ and $t$, so
equation \eqref{E1} and Proposition \ref{7sumofpowers} tell us that
deg($\aaa$) = 2 deg($f$)+2 in $m$, and 2 deg($f$)+1 in $r$.
\end{proof}\hfill

\section{The Poincare-Birkhoff-Witt theorem for $H_f$}

The proof of the PBW theorem below, builds on Section \ref{S11} above. We
first remark, though, that the PBW theorem (and hence the analysis in
Section \ref{S11}) can all be proved using the Diamond Lemma (cf.
\cite{Be}). This (was suggested by W.L. Gan to the author, and) is done
in detail in future work, with W.L. Gan and N. Guay, in \cite{GGK}, for a
similar associative algebra - namely, the $q$-analog of $H_f$.\\

We now show the PBW theorem for $H_f$. If $\dn = 0$, then $H_f$ is the
universal enveloping algebra of a five-dimensional Lie algebra, so the
PBW theorem holds. If not, then to show the PBW theorem, we need the
following key lemma.

\begin{lemma}\label{L4}
Given $s \in \nn$, there is a finite subset $T \subset k$ so that if $r
\notin T \cup \nn$, then $X^s v_{r-s} = X^s p_s(Y,F) \neq 0$ in $Z(r)$.
\end{lemma}

\noindent (Note that since char $k$ = 0, hence $\Z \hookrightarrow k$,
and therefore $\nn \cup\ \{$a finite set\} $\neq k$.)\\

\begin{proof}
If $r \notin \nn$, then repeatedly applying $(R_t)$ yields
$$X^{s-1} v_{r-s} = X^{s-1} p_s(Y,F) = [(r-s+2)(r-s+3) \dots r]^{-1}
(-1)^{s-1} [\aaa_{r,s-1} \aaa_{r,s-2} \dots \aaa_{r,3}]\ v_{r-1}$$

The first product of terms is nonzero if we take $r \notin \nn$, so
denote it by $d_0 \neq 0$. Also, $X v_{r-1} = XY v_r = -\dn v_r = -c_{0r}
v_r$. Therefore
$$X^s v_{r-s} = X^s p_s(Y,F) = (-1)^s\bigg( d_0 c_{0r}
\prod_{j=3}^{s-1} \aaa_{r,j} \bigg) v_r$$

Clearly each term in the product is a polynomial - but this time in $r$
(by Corollary \ref{C11}), as is $c_{0,r}$ (by definition). Therefore let
us take $T$ to be the set of roots of {\it all} these polynomials in $k$.
Clearly, if $r \notin \nn \cup T$, then the right hand side does not
vanish in $Z(r) \neq 0$, and hence we are done.
\end{proof}\hfill

We prove two claims, and then the PBW theorem. As above, we take $(N_+)$
to be the left ideal generated by $N_+ = kX \oplus kE$. But first, we
observe that $B_- = \mf{h} \bigoplus N_- = kH \oplus kY \oplus kF$ is a
Lie algebra, so we know the PBW theorem for it. Consequently, the
multiplication map: $k[Y,F] \otimes_k k[H] \to \mf{U}(B_-)$ is an
isomorphism.

\begin{prop}
$k[Y,F] k[H] \bigcap (N_+) = 0$.
\end{prop}

\begin{proof}
Suppose $\exists 0 \neq b \in k[Y,F] k[H] \cap (N_+)$. Now, $b=0$ in
every Verma module $Z(r)$, so $b_+ b$ is also zero, for every $b_+ \in
\mf{U}(N_+) = k[X,E]$.

But we will now produce $b_+$ and $r$ so that $0 \neq b_+ b \in k^\times
\cdot \bar{1}$ in $Z(r)$, thus producing a contradiction. Suppose $b_-$
is of the form $\sum_{i,j} Y^i F^j b_{ij}(H) \in k[Y,F]k[H]$. Firstly, we
may assume w.l.o.g. that $b_-$ is a weight vector for $H$, because if
not, then we take the lowest weight component to $k^\times \cdot 1$, and
then the other components automatically are killed.\\

So suppose $b_- = \sum_{j=0}^l F^j Y^{n-2j} b_j(H)$. Let $l'$ be the
largest number so that $b_{l'}$ is nonzero. W.l.o.g. $b_l \neq 0$ (i.e.
$l'=l$), so $b_l$ has a finite set of roots $S$. Also, given $l$, the
above lemma says there exists a finite set $T$ so that if $r \notin \nn
\cup T$, then $X^{n-2l}v_{r-(n-2l)} \in k^\times v_r = k^\times \cdot
\bar{1}$.

So fix $r \notin \nn \cup T \cup S$. Then $b_- = \sum_{j=0}^l Y^{n-2j}
F^j b_j(r)$, and $b_l(r) \neq 0$. We now write $b_-$ as a linear
combination
$$b_- = a_0 v_{r-n} + a_2 F v_{r-n+2} + \dots + a_{2l} F^l v_{r-n+2l}$$
where $a_{2l} = b_l(r) \neq 0$, because $v_{r-n} = p_n(Y,F)v_r$, and the
$p_n$'s are monic in $Y$.\\

Since the $v_t$'s are \spl-maximal, hence by \spl-theory, $E^l$ kills all
summands but the last one. And since $r \notin \nn \cup T \cup S$, hence
again by $\spl$-theory (cf. \cite[$\S 7 \S$]{H}), $E^l b_- = E^l (a_{2l}
F^l v_{r-n+2l}) = c_0 v_{r-n+2l}$ for some nonzero scalar $c_0$. But then
$X^{n-2l}(E^l b_-) = c_0 X^{n-2l} v_{r-(n-2l)}$, and this is nonzero by
the above lemma. Hence we have produced $b_+$ so that $b_+ b \neq 0$ in
$Z(r)$. This is a contradiction to the first paragraph in this proof, and
hence we are done.
\end{proof}\hfill

\begin{cor}
$Z(r) \cong k[Y,F]\ \forall r \in k$.
\end{cor}

\begin{proof}
Suppose not. Then there is a relation, say of the form $b_- \in k[Y,F]
\cap (N_+, (H - r \cdot 1 ))$. Since the multiplication map: $k[Y,F]
\otimes_k k[H] \otimes_k k[X,E] \to H_f$ is onto, hence say $b_- = n_+ +
p$, where $n_+ \in (N_+)$, and $p \in k[Y,F] k[H] \setminus k[Y,F]$.
Clearly, then, $n_+ = b_- - p \in k[Y,F] k[H] \cap (N_+) = 0$.

Further, $p$ is of the form $p = \sum_i b_{-i} p_i(H - r \cdot 1)$, where
each $p_i$ is a polynomial with no constant term, and the $b_{-i}$'s are
linearly independent in $k[Y,F]$. Since we know the PBW theorem for the
{\it Lie algebra} $B_-$, hence $k[Y,F] \otimes_k k[H] \cong k[Y,F]k[H]$.
Thus $p_i = 0\ \forall i$, so $p=0$, whence $b_- = 0$ as required.\\
\end{proof}

\noindent Finally, we have

\begin{theorem}\label{9pbw1}
The PBW theorem holds, i.e. $\{ F^a Y^b H^c X^d E^e : a,b,c,d,e \geq 0
\}$ is a $k$-basis for $H_f$.\\
\end{theorem}

\begin{proof}
Suppose not. Then there is a relation of the form $a = \sum_{i=1}^l b_i
X^{d_i} E^{e_i} = 0$, where $b_i \in k[Y,F] k[H]$ for each $i$.\\

We first find $b_- \in k[Y,F]$ on which {\it exactly one} of the
$X^{d_i}E^{e_i}$'s acts nontrivially. Choose the least $e$, and among all
$d_i$'s, choose the least $d$, for which $X^d E^e$ has nonzero
coefficient. By the above lemma, there exists a finite set $T$ so that
$X^d v_{r-d} \neq 0$ in $Z(r)$ if $r \notin \nn \cup T$.

Let us now look at $v= F^e v_{r-d} \in k[Y,F]$. Clearly, for $(d',e')
\neq (d,e)$, either $e' > e$ (in which case $(X^{d'} E^{e'}) (F^e
v_{r-d}) = c_0 (X^{d'} E^{e'-e-1}) Ev_{r-d} = 0$), or $e' = e$ and $d' >
d$ (in which case $(X^{d'} E^e) (F^e v_{r-d}) = c_0 X^{d'} v_{r-d} = c_0'
X^{d'-d-1} Xv_r = 0$), for some nonzero $c_0, c_0' \in k$. Thus we see
that only $X^d E^e$ acts nontrivially on $v \in Z(r)$, because $(X^d
E^e)(F^e v_{r-d}) = c_0 X^d v_{r-d} = c_0' v_r$ for $c_0,c_0' \in
k^\times$, from above. Thus we have found such a $b_- \in k[Y,F]$.\\

Returning to the PBW theorem, recall that we had a linear combination
that was zero: $a = \sum_{i=1}^l b_i X^{d_i} E^{e_i} = 0$, and w.l.o.g.
we assume the special $(d_i,e_i)$ (as above) corresponds to $i=l$. Now
suppose that $b_l = \sum_j b_{-j} p_j(H)$, where $b_{-j}$ are linearly
independent in $k[Y,F]$, and $p_j$ are nonzero polynomials. Then $\Pi p_j
= p \neq 0$, and $k \setminus (\nn \cup T)$ is infinite, so choose any $r
\notin (\nn \cup T)$, such that $p(r) \neq 0$. Therefore $p_j(r) \neq 0\
\forall j$.\\

Finally, we have $a = 0$, so $0 = a \cdot b_-$ (where $r$ is chosen
above) $= c_r b_l$ for some nonzero scalar $c_r$ (note that we are
working in $Z(r)$ here). Therefore $b_l$ is zero in $Z(r)$, whence
$\sum_j p_j(r) b_{-j} = 0$. But the $b_{-j}$'s are linearly independent
in $Z(r) \cong k[Y,F]$ (from above), and $p_j(r) \neq 0\ \forall j$ (by
choice of $r$). This is a contradiction, hence such a relation $a=0$
cannot occur in the first place.
\end{proof}\hfill

\section{Necessary condition for $Z(t) \hookrightarrow Z(r)$}

\noindent The main result is

\begin{theorem}\label{Tnec}
\hfill
\begin{enumerate}
\item If $Z(r)$ has a maximal vector of weight $r-n=t$, then (it is
unique upto scalars, and) $\aaa_{r,r-t+1} = 0$.
\item (\emph{Verma's Theorem}, cf. \cite{Ver}, \cite[(7.6.6)]{Dix})
$\HHom(Z(r'),Z(r)) = 0$ or $k$ for general $r,r' \in k$. All nonzero
homomorphisms are injective.
\end{enumerate}
\end{theorem}

The first part of Verma's theorem is easy to show given the previous
part, and the second part follows from Corollary \ref{C5}. For the first
part of the theorem, we need some preliminaries.

\noindent {\bf Definition :} Given $T \in H_f$, denote by $W(r,n,T)$ the
set of solutions to $Tv = 0$ in $Z(r)_{r-n}$.\\

\begin{prop}\label{11general1}
For all $n \in \nn$ and $r \in k$, we have
\begin{enumerate}
\item $\dim_k(W(r,n,X)) \leq 1$; it equals 1 if $n$ is even.
\item $1 \leq \dim_k(W(r,n,E)) \leq 2$ if $r+1 \in \nn$ and $r+1 \leq n
\leq 2r+2$; it equals 1 otherwise.
\end{enumerate}
\end{prop}

\begin{proof} Both the proofs are similar, so we show (1) now. We know
$Z(r)_{r-n}$ is spanned by $Y^n, FY^{n-2},\dots$. Now, we claim that if
$Xv = 0$ for nonzero $v \in Z(r)_{r-n}$, then the contribution of $Y^n$
to $v$ is nonzero (i.e. $v = a_0 Y^n + a_1 F Y^{n-2} + \dots$, where $a_0
\neq 0$).

Well, suppose $v = \sum_{i \geq s} a_i F^i Y^{n-2i}$ for some $s \geq 0$,
where $a_s \neq 0$. From Lemma \ref{Lbasis1}, we see that $Xv = -sa_s
F^{s-1} Y^{n-2s+1} +$ terms of lower degree in $Y$. Since $a_s \neq 0$,
hence $s=0$ as required.

Thus, every $0 \neq v \in W(r,n,X)$ is of the form $v = c Y^n +$ lower
order terms. Now suppose we have two such $0 \neq v_i = c_i Y^n + l.o.t.
\in W(r,n,X)$ (i.e. for $i=1,2$). Then $c_2 v_1 - c_1 v_2$ is also in
$W(r,n,X)$, but without any $Y^n$ term. Hence it is zero from above, so
that $v_2 \in k \cdot v_1$, as required.

Finally, we need to show that if $n$ is even, then such a $v$ exists.
Recall the Kostant function $p$. Now observe that $p(-2n) = p(-2n+1) + 1\
\forall n$ (because we have the sets $\{ F^0 Y^{2n}, \dots, F^n Y^0 \}$
and $\{ F^0 Y^{2n-1}, \dots, F^{n-1} Y \}$). Thus, $X : Z(r)_{r-2n} \to
Z(r)_{r-2n+1}$ is a map from one space to another of lesser dimension.
Hence it has nontrivial kernel, as required.
\end{proof}\hfill

\begin{remark}\hfill
\begin{enumerate}
\item This makes the relation $Xv_t \in k v_{t+1}$ easier to understand:
$E(Xv_t) = X(Ev_t) = 0$, so $Xv_t$ is in $W(r,r-t-1,E)$.
\item The above result holds for any $Z(r) \to V \to 0$. In any such $V$,
any maximal vector of a given weight $r'$ (if it exists)is unique upto
scalars.\\
\end{enumerate}
\end{remark}

\begin{prop}\label{12general2}
We work again in the Verma module $Z(r)$ for any $r \in k$.
\begin{enumerate}
\item $\dn$ acts on $F^m Y^n$ by $\dn F^m Y^n = F^m(c_{0,r-n} Y^n +
l.o.t.) \in Z(r)_{r-n-2m}$.
\item If $v \in Z(r)_{r-n}$ satisfies $Xv = 0$, then upto scalars we have
$$v = Y^n - FY^{n-2} \sum_{l=0}^{n-1} c_{0,r-l} + l.o.t.$$
\item If $v \in Z(r)_{r-n}$ satisfies $Ev = 0$, then upto scalars, $v$ is
one of the following:\\
\hspace*{3.6ex} (a) $v = F^{j+1} v_j$, where $-1 \leq j \leq r,\ r+1 \in
\nn$, and $r+1 \leq n \leq 2(r+1)$\\
OR (b) $\displaystyle v = (r+2-n)Y^n + FY^{n-2} \sum_{m=0}^{n-2} (n-1-m)
c_{0,r-m} + l.o.t.$
\end{enumerate}
\end{prop}

\begin{remark}\hfill
\begin{enumerate}
\item Here, $l.o.t.$ denotes monomials of lower order in $Y$.
\item Thus, a necessary condition for $Z(r)$ not to be simple (for
general $r \notin k$) is that $\aaa_{r,r-t+1} = 0$ for some $t \in r -
\N$. Further, if $r \notin \nn$, then Corollary \ref{C6} says that this
condition is also sufficient, i.e. the converse to (4) holds as well, if
the maximal vector $v_t$ is nonzero.\\
\end{enumerate}
\end{remark}

\begin{proof}\hfill
\begin{enumerate}
\item W.l.o.g. $m=0$, because $\dd$ (and hence $\dn$) commutes with $F$.
We now proceed by induction on $n$. For $n=0,\ v_r$ is maximal, hence
(e.g. cf. Corollary \ref{C11}) $\dn v_r = c_{0r} v_r$. Further, $\dn = 1
+ f(\dd)$ and hence $\dn \in \text{End}_k(Z(r)_t)$ for any $t \in r -
\nn$.

Thus, $\dn Y^n$ is a linear combination of $Y^n,\ FY^{n-2}$, and lower
order terms in $Y$. Now, $4\dd = 2FE + (H^2+2H)/2$, so $4\dd Y^n = 2FE
Y^n + [(H^2+2H)/2] Y^n$. Of course, $EY^n$ is a linear combination of
$Y^{n-2-i} \dn Y^i$ from above, and $\dn Y^i$ is a linear combination of
lower order terms, by induction. So $EY^n$ and hence $2FEY^n$ are l.o.t.
in $Y$.

Thus, $\dd Y^n = [(r-n)^2 + 2(r-n)]Y^n / 8\ + l.o.t. = c_{r-n} Y^n +
l.o.t.$ (because $H$ acts on $Z(r)_{r-n}$ by $r-n$). Also, we have $\dd
(l.o.t.) = l.o.t.$ by the induction hypothesis, so $\dd^2 Y^n = c_{r-n}^2
Y^n + l.o.t.$, and so on. Hence $\dn Y^n = (1 + f(\dd)) Y^n = (1 +
f(c_{r-n})) Y^n + l.o.t. = c_{0,r-n} Y^n + l.o.t.$ as required.\\

\item From Lemma \ref{Lbasis1}, $\displaystyle XY^n = - \sum_{l=0}^{n-1}
Y^{n-1-l} \dn Y^l = - Y^{n-1} \sum_{l=0}^{n-1} c_{0,r-l} + l.o.t.$ by
what we just proved. Similarly, $XFY^{n-2} = -Y^{n-1} + l.o.t.$, and
hence if $Xv = 0$, then $v$ is monic in $Y$, and it {\it must} look like
$\displaystyle v = Y^n - FY^{n-2} \sum_{l=0}^{n-1} c_{0,r-l} + l.o.t.$,
in order that the two highest degree (in $Y$) terms vanish.\\

\item The argument is the same as the one just above; the coefficients
are slightly different.
\end{enumerate}
\end{proof}

\begin{proof}[Proof of Theorem \ref{Tnec}]
If $v = Y^n + l.o.t. \in Z(r)_{r-n}$ is maximal, then so is $(r-n+2)v$,
and then both conditions (the ones in (2) and (3)(b) above) must be
satisfied, whence the coefficient of $FY^{n-2}$ is the same in both the
forms. Therefore we have
$$-(r-n+2)\sum_{l=0}^{n-1} c_{0,r-l} = \sum_{m=0}^{n-2} (n-1-m)c_{0,r-m}
= \sum_{l=0}^{n-1} (n-1-l)c_{0,r-l}$$
because for $l=n-1$ the summand on the RHS vanishes. Simplifying this,we
get\\
$\displaystyle \sum_{l=0}^{n-1} [(r-n+2)+(n-1-l)] c_{0,r-l} = 0$, which
by definition means $\aaa_{r,n+1} = \aaa_{r,r-t+1} = 0$ as required.
\end{proof}\hfill

Suppose $\dn \neq 0$. Given $r \in k$, let $r_0$ be the maximal $t \in r
+ \nn$, such that $t=r$ is a root of $\aaa_{r_0,r_0-t+1}$ (this exists
because $\aaa_{rm}$ is a polynomial, as in Corollary \ref{C11}). Define
the set $S(r)$ to be the set of roots $t$ of $\aaa_{r_0,r_0-t+1}$, that
are in $r_0 - \nn$.

We claim that if $\aaa_{t,t-t'+1} = 0$, then $t \in S(r)$ iff $t' \in
S(r)$. (Thus, $S(r)$ is the transitive (and symmetric) closure of
$\{r\}$, under the relation of ``being a root of $\aaa_{t,m}$".) This
follows from equation \eqref{eqn1} (mentioned in the proof of Corollary
\ref{C9}).\\

\begin{lemma}\label{L2.2}
Suppose $\dn \neq 0$.
\begin{enumerate}
\item For any $r \in k$, the set $S(r)$ is finite, of size at most $2\
deg(f) + 2$.
\item The sets $S(r)$ partition $k$.
\end{enumerate}
\end{lemma}

\begin{proof}
The first part follows from Corollary \ref{C11}, and the second part from
equation \eqref{eqn1}.
\end{proof}

\noindent {\bf Warning}: The set $S(r)$ need {\it not} serve the role of
the $S(\la)$'s of the first part (of this paper), but might split into a
disjoint union of sets $S(\la)$. As we shall see later, in most cases the
$S(r)$'s do serve as $S(\la)$'s, though.\\

\section{Finite dimensional simple $H_f$-modules}

Suppose $V = V(r)$ is finite-dimensional and simple. Then $r \in \nn$,
and $V = \oplus V_C(n)$, as mentioned earlier (or cf. \cite[$\S 7.8
\S$]{We}). (Here, $0 \leq n \leq r$ for each summand.) Thus any nonzero
\spl-maximal weight vector in $V(r)$ has non-negative weight. In
particular, $v_{-1}=0$ in $V(r)$.

The highest weight space has dim$_k(V_r) = 1$, so $[V(r) : V_C(r)] = 1$.
Let us use equations \eqref{Rt}, \eqref{St} now. We know $v_{-1} = 0$ in
$V(r)$, so let $s$ be the largest integer in $\nn$ such that $v_{s-1} =
0$ but $v_s$ is nonzero in $V(r)$. Thus, $v_t \neq 0$ if $s \leq t \leq
r$ by Corollary \ref{C7}. Also, by Corollary \ref{C6}, we have
$\aaa_{r,r-s+2}$ (and hence $d_{s-2}$ if $s > 0)\ = 0$ (but $d_t \neq 0\
\forall t \in s - 1 + \nn$). Thus $Yv_s = -d_{s-1} Fv_{s+1}$ etc. Now,
the equations \eqref{Rt} and \eqref{St} show us that the subspace
$\bigoplus_{i=s}^r V_C(i)$ is an $H_f$-submodule of $V(r)$. Since $V(r)$
is simple, they are equal, and we have just proved\\

\begin{theorem}\label{4finsimple}
If $V = V(r)$ is finite-dimensional, then $r \in \nn$ and $\exists s \leq
r \in \nn$ so that $V = \bigoplus_{i=s}^r V_C(i)$. Also, $\aaa_{r,r-s+2}
= v_{s-1} = 0$ and $\Pi(V) = \{ \pm r ,\ \pm (r-1),\ \dots,\ \pm s \}$ is
$W$-stable. Conversely, if $\exists 0 \leq s \leq r$ so that
$\aaa_{r,r-s+2} = 0$, but $d_t \neq 0\ \forall s-2 < t < r-1$, then $V =
\bigoplus_{i=s}^r V_C(i)$, where $V_C(i)$ is a simple \spl-module with
\spl-maximal vector $v_i$.
\end{theorem}

\begin{remark}
The module structure is completely determined by the relations
\eqref{Rt}, \eqref{St}, and \spl-theory.\\
\end{remark}

\noindent (The Weyl group $W$ acts on $\Pi(V)$ (and $V$) by permuting $\{
\mu,\ -\mu \}$ (and $\{ V_\mu,\ V_{-\mu} \}$), as seen in the next
section.) We say an ideal $I$ of $H_f$ is {\it primitive} if $H_f / I$ is
a simple $H_f$-module. Define $J(r)$ to be the annihilator of $V(r) =
V(r,s)$ in $H_f$ (we still have $r \in \nn$, of course), and let $Y(r) =
\Rad(Z(r))$.\\

\begin{prop}\label{5maxsub}
$J(r)$ is generated by $\{ F^{j+1} p_{r-j}(Y,F) = F^{j+1}v_j : s \leq j
\leq r \}$ along with $p_{r-s+1}(Y,F) = v_{s-1}, N_+,$ and $(H-r \cdot
1)$. Further, if $j \in \nn$ then we have $X F^{j+1}v_j = -(j+1)F^j
v_{j-1}$.
\end{prop}

\begin{proof}
Observe that $J(r)$ definitely contains all these terms because these
relations vanish in $V(r,s)$ (where $1 = v_r$). So let these relations
generate the (left) ideal $I$. (Thus $H_f / I \twoheadrightarrow V(r,s) =
H_f / J(r) = Z(r) / Y(r)$.) Since $p_{r-s+1} \in I$, we see that every
element in $H_f / I$ is of the form $Y^j F^k$ where $j \leq r-s$. Then
the other relations tell us that $0 \leq k \leq r-j$. Thus $\dim_k (H_f /
I) \leq (r+1) + r + \dots + (s+1)$, and we can easily verify (using
Theorem \ref{4finsimple}) that this {\it is} $\dim_k(V(r,s))$. Hence we
are done.

For the second part, we calculate: $[F^n,X] = n F^{n-1} Y$. Then the
rest is (also) calculation.\\
\end{proof}

\section{Characters, and an automorphism}\label{sec19}

Recall that we have already defined the group ring $\Z[\Z]$, Kostant and
Weyl's functions, and the formal character earlier. Since we know that
$Z(r) \cong k[Y,F]$ as $\mf{U}(N_-)$-modules, hence $p(-n) = 1+ \lfloor
n/2 \rfloor$ for $n \in \nn$, if we identify $\eta$ in $\Phi_f^+$ with 1,
and hence $2 \eta$ with 2. The {\it Weyl function} $q$ is just $(e(1) -
e(-1))(e(1/2) - e(-1/2))$. Also define

\begin{center}
$\omega(r + \delta,s + \delta') = \biggl[ \sum_{\sigma \in W}sn(\sigma)
e(\sigma(\frac{r+s+2}{2})) \biggr] \biggl[ \sum_{\sigma \in W}sn(\sigma)
e(\sigma(\frac{r-s+1}{2})) \biggr]$\\
\end{center}

\noindent where $\delta = 3/2 = 3 \eta / 2 = \frac{1}{2} \sum_{\aaa \in
\Phi_f^+} \aaa$, and $\delta' = \eta / 2$. Thus we have
$\omega(r+\delta,\ s+\delta') = e(r+\delta) - e(s+\delta') -
e(-s-\delta') + e(-r-\delta)$.\\

\begin{lemma}\label{L3} \hfill
\begin{enumerate}
\item $ch_{Z(\la)} = e(\la)(1+t)(1 + 2t^2 + 3t^4 + \dots) = p *
\epsilon_\la$, where $t = e(-1)$.
\item $q = \omega(\delta,\ \delta')$.
\end{enumerate}
\end{lemma}

\noindent The proof is a matter of easy calculation.\\

\noindent {\bf Digression on $W'$:}\\
We now discuss the action of a different group $W' = (\Z / 2 \Z)^2$ on
the roots. We work here with $\Z_2$, the ring of {\it dyadic fractions}
$\{ a / 2^b:\ a,b \in \Z \}$. First of all, $W' = \{ 1, \ \sigma_1,\
\sigma_2,\ \sigma_1 \sigma_2 = \sigma_2 \sigma_1 = -1 \}$. Further, it
acts on $M = \halfz \times \halfz$ by invertible linear maps, i.e. $W'
\subset GL_2(\halfz)$.

To compute the explicit action, let $e_1 = (1,0),\ e_2 = (0,1)$ be a
$\halfz$-basis for the free $\halfz$-module $M$ of rank 2. Then
$\sigma_i(e_j) = (-1)^{\delta_{ij}} e_j$ where $i,j \in \{ 1,2 \}$.
Further, there is a {\it sign} homomorphism $\sn : W' \to \{ \pm 1 \}$,
given by $\sn(\sigma) = \det(\sigma) = (-1)^{l(\sigma)},\ l$ being the
{\it length}. Thus the $\sigma_i$'s are transpositions, or more
accurately, reflections, and sn($\sigma_1 \sigma_2)$ = 1, because
$\sigma_1 \sigma_2 = (-1) \cdot \id$ on all of $M$.\\

We now define a map $\varphi : \halfz \times \halfz \to \halfz$, given by
$\varphi(m,n) = m+2n$. This corresponds to identifying the first
coordinate with the coefficient of $\eta$ = root of $X$, and the second
with the coefficient of $2 \eta$ = root of $E$. Thus, the half sum of the
roots would be $\delta = \frac{1}{2}\varphi(1,1) = \varphi(\half,\half)$.
Similarly, $\varphi(-1,1) = 2 \delta'$.

Given $(m,n) \in \halfz \times \halfz$, we draw a ``square" of its orbit
under $W'$. In what follows below, we ``cut" off a side of the square and
expand out the sides in one line. For instance, we have the following map
(essentially, we {\it want} this to hold, in order to write formulae for
$V(r,s)$ analogous to the \spl-case):
$$\start{\half}{\half} \diag{\delta}{\delta'}{-\delta}{-\delta'}$$

Identifying $\delta, \eta, \delta'$ etc. with numbers in $\Z$, we can
write
\begin{equation}\label{E2}
\start{\half}{\half} \diag{3/2}{1/2}{-3/2}{-1/2}
\end{equation}
\hfill

The orbit of $W'$ - or more precisely, $\varphi \circ W'$ - on the roots
$\eta = e_1$ and $2 \eta = e_2$, is given by
\begin{equation}\label{E3}
\start{1}{0} \diag{1}{-1}{-1}{1}
\end{equation}
\begin{equation}\label{E4}
\start{0}{1} \diag{2 \eta}{2 \eta}{-2 \eta}{-2 \eta}
\end{equation}

And then we see that \eqref{E3} + \eqref{E4} = \eqref{E2} + \eqref{E2},
which {\it should} hold, because we defined $\delta$ as the half sum of
positive roots - and which {\it does} hold, because the actions of
$\varphi$ and $\sigma \in W'$ are all linear.\\

Finally, if $V = V(r,s)$ is a simple $H_f$-module, then we also have
\begin{equation}\label{E5}
\start{\frac{r-s+1}{2}}{\frac{r+s+2}{4}} \diag{r + \delta}{s +
\delta'}{-r - \delta}{-s - \delta'}
\end{equation}

Note that equation \eqref{E2} is a special case of this last equation
\eqref{E5}, if we take $r=s=0$. Now denote by $\psi$ the endomorphism of
$\halfz \times \halfz$, sending $(r,s)$ to
$(\frac{r-s+1}{2},\frac{r+s+2}{4})$. We can now use this to write the
standard character formulae.\\

\noindent {\bf Back to characters :}\\
We see now that $\displaystyle \omega(r + \delta,s + \delta') =
\sum_{\sigma \in W'} sn(\sigma) e(\varphi \sigma \psi(r,s))$ and
$\omega(\delta,\delta') = \sum_{\sigma \in W'} sn(\sigma) e(\varphi
\sigma \psi(0,0))$.\\

Let us now look at $ch_{V(r,s)} = ch_{r,s}$, say, where $V(r,s)$ is {\it
simple}. Theorem \ref{4finsimple} says that $ch_{r,s} = \sum_{i=s}^r
ch(V_C(i))$, and so we have (exactly as in \spl-theory)\\

\begin{theorem}\label{6wcf}
Say $V = V(r,s)$ is a simple $H_f$-module. Then we have\\
\begin{enumerate}
\item (Weyl's Character Formula)
$$\omega (\delta, \delta') * ch_{r,s} = \omega (r + \delta, s +
\delta'),\ or\ ch_{r,s} = \frac{\sum_{\sigma \in W'} sn(\sigma) e(\varphi
\sigma \psi(r,s))}{\sum_{\sigma \in W'} sn(\sigma) e(\varphi \sigma
\psi(0,0))}$$

\item (Alternate version of the Weyl Character Formula)
$$e(\delta) ch_{r,s} = \omega(r+\delta, s+\delta') * ch_{Z(0)} =
\sum_{\sigma \in W'} sn(\sigma) ch_{Z(\varphi \sigma \psi (r,s))}$$
%

\item (Kostant's Multiplicity Formula) Say $m_r (t) = dim(V(r,s)_t).$
Then
$$m_r (t) = (p * \epsilon_{-\delta} * \omega(r+\delta, s+\delta'))(t) =
\sum_{\sigma \in W'} sn(\sigma) p(t + \delta - \varphi \sigma \psi
(r,s))$$

\item (Weyl's Dimension Formula)
$$\deg(r,s) \ (\ \overset{ \text{def} }{=} \dim V(r,s)) =
\lim\limits_{e(1) \to 1} ch_{r,s} = \frac{(r+s+2)(r-s+1)}{2} =
\frac{\psi_1(r,s)\psi_2(r,s)}{\psi_1(0,0)\psi_2(0,0)}$$
where $\psi = (\psi_1,\psi_2)$.\\
\end{enumerate}
\end{theorem}

\section{Standard cyclic modules for $r \not\in \nn$}

\begin{stand}
For the rest of this paper, we assume that $\dn \neq 0$.
\end{stand}

We now examine the structure of standard cyclic modules $Z(r) \to V \to
0$, for various $r \in k$. The easier choice is $r \notin \nn$. Theorem
\ref{3basics} says that the equations \eqref{Rt}, \eqref{St} are valid
for all $t \in r - 2 - \nn$, so we can define the \spl-maximal vectors
$v_t$ for all $t$. Theorem \ref{3basics} tells us that these span all the
\spl-maximal vectors.\\

Hence the only maximal vectors in $V$ are those $v_t$'s for which
$\aaa_{r,r-t+1} = d_{t-1}=0$. (Thus there are finitely many maximal
vectors.) Now say $W$ is a submodule of highest weight $t$ for some such
$t$. We claim that $W = Z(t)$. Suppose not, i.e. say $W$ contains a
vector of the form $a_1 F^{i_1} v_r + \dots + a_m F^{i_m} v_{t+1}$ (in
addition to $Z(t)$). Repeatedly applying $E$, we conclude that $W$
contains a vector of weight higher than $t$, a contradiction. (We use
similar arguments in $\S \ref{sec21} \S$ below.) Thus there are finitely
many submodules, and $V$ has a finite composition series, given by the
{\it distinct} roots of $\aaa_{r,m}$ that are in $r - \nn$.

\begin{theorem}\label{8rnotinnn}
Suppose $r \notin \nn$, and $Z(r) \to V \to 0$.
\begin{enumerate}
\item The only submodules of $V$ are $H_f v_t = \mf{U}(N_-) v_t$, where
$t = r-m+1$ is a root of $\aaa_{rm}$, i.e. $\aaa_{r,r-t+1} = d_{t-1} =
0$. These are only finitely many.
\item $V$ has a unique composition series with length at most $\emph{deg
} \aaa_{rt} = 2(\emph{deg}(f)+1)$.
\item The composition factors are isomorphic to $Z(t_i) / Y(t_i) =
V(t_i)$, one for each root $t_i \in r - \nn$ and nonzero maximal vector
$v_{t_i}$.
\item Given $r' \in k,\ \HHom(Z(r'),Z(r)) \neq 0$ iff $r' = t_i$ for some
$i$.
\item The primitive ideal here is generated by $v_{t_1} =
p_{r-t_1} (Y,F)$ (for the ``largest" such $t_1$).\\
\end{enumerate}
\end{theorem}

\section{Standard cyclic modules for $r \in \nn$}\label{sec21}

We now consider the case when $r \in \nn$. Let $r = t_0 > t_1 > \dots >
t_k \geq -1$ be all the distinct integers so that $v_{t_j}$ is a maximal
vector in $Z(r)$ (i.e. all the distinct roots ($\geq -1$) of
$\aaa_{r,r-t+1}$). We define the $H_f$-submodule $Y(t_i,t_j)$ to be the
\usl-submodule generated by $\{ F^{m+1} v_m : t_i \leq m \leq t_j \}$,
and $Z(t_i)$. Clearly, we have $t_i \leq t_j$, or $i \geq j$, and we also
have the obvious inclusions $Y(t_i,t_j) \subset Y(t_{i'},t_{j'})$ iff
$t_i \leq t_{i'}$ and $t_j \leq t_{j'}$.

Now if $V(r) = V(r,s)$ is simple, then $r=t_0, s=t_1$. Also, we clearly
have $Z(t_i) = Y(t_i,t_i)$ and $Y(t_i) = Y(t_{i+1},t_i)$ is the maximal
submodule of $Z(t_i)$. We now classify some submodules of $Z(r) = Z(t_0)
= Y(t_0,t_0)$, and show that $Z(r)$ has finite length.\\

\begin{prop}\label{P3}
$Y(r) = Y(t_1,t_0)$, and every submodule of $Z(r)$ is either of the form
$Y(t_l,t_s)$ (for some $k \geq l \geq s \geq 0$), or all its weights are
(strictly) below $t_k$.
\end{prop}

\noindent {\bf Proof : (a)} Suppose $V$ is a submodule. We first show
that if $F^{j+1}v_j \in V$ (for some $j \geq -1$), then $V$ is of the
form $Y(t_l,t_s)$ for some $k \geq l \geq s \geq 0$.\\

Suppose $F^{j+1}v_j \in V$. Then $V$ also contains $XF^{j+1} v_j =
-(j+1)F^j v_{j-1}$ (by Proposition \ref{5maxsub}), and repeatedly
applying $X$, we conclude that $v_{-1} \in V$. Keep on applying $X$, to
get that $v_0, v_1$, and so on are in $V$, until $v_{t_k} \in V$, because
this is the first point where we cannot get further ahead (because
$d_{t_k- 1} = \aaa_{r,r-t_k +1} = 0$). Thus, if $v'$ is a weight vector
of highest possible weight $x$ in $V$, then $x \geq t_k \geq -1$. Also,
$Ev' = 0$, meaning that $v' = v_x$ upto scalar, from part (4) of Theorem
\ref{3basics}. Next, $Xv' = Xv_x = 0$, so $d_x = 0$, meaning that $x =
t_l$ for some $l$ (by Corollary \ref{C6}).\\

Thus, if $F^{j+1}v_j \in V$ for some $j$, then $V$ contains $Z(t_l)$ as
well as the \usl-span of $F^{j+1}v_j$'s, say for $0 \leq j \leq m (\leq
r)$ ($m$ maximal). Again, if $F^{m+1}v_m \in V$, then $XF^{m+1}v_m =
-(m+1)F^mv_{m-1} \in V$, and as above, $YF^{m+1}v_m = F^{m+1}(v_{m-1} -
d_{m-1}Fv_{m+1}) \in V$. But now, $d_{m-1} = 0$ iff $v_m$ is maximal (by
Corollary \ref{C6}). Thus if $v_m$ is not maximal then $F^{m+2}v_{m+1}
\in V$ as well. But $m$ was chosen to be maximal; hence $v_m$ has to be
maximal, and $m = t_s$ for some $s$. Thus we conclude that $Y(t_l,t_s)
\subset V$.\\

If this inclusion is proper, then $V$ contains a linear combination of
terms of the form $F^{j+1+m} v_j\ (m \geq 0, j > t_s)$ and $F^m v_i\ (0
\leq m \leq i, i > t_l)$. Since all $F^{j+1}v_j$'s and $v_i$'s are
\spl-maximal, hence repeatedly applying $E$ gives that a linear
combination of $F^{j+1}v_j$'s and $v_i$'s is in $V$. Now we use the
$H$-action to separate all these terms, and we conclude that $V$ contains
a term of the form $F^{j+1}v_j$ for $j>t_s$, or $v_i$ for $i>t_l$. This
contradicts the maximality of $t_s,t_l$, hence $V = Y(t_l,t_s)$ as
claimed.\\

\noindent {\bf (b)} Now, if $V$ contains no vector of the form
$F^{j+1}v_j$ (for $-1 \leq j \leq r$), then we claim that $V$ has weight
vectors with weights only below $t_k$. For if not, then $V$ contains a
vector in the \usl-span of higher weight vectors $v_t (t_k \leq t \leq
r)$, which would mean it would contain $F^i v_j$ for some $i,j$ (by
similar application of $E,H$ as above), and multiplying by a suitable
power of $F$ gives us that $F^{j+1}v_j \in V$ for some $j$. This is
false. \qed\\

In general, we know that either $Z(t_k)$ is simple, or $Y(t_k)$ has a
maximal vector of highest possible weight $t$, say, which is $\leq -2$.
We now find all submodules of $Z(t_k)$, or equivalently, of $Y(t_k)$. (Of
course, if $t_k = -1$ then we are already done, because $Z(-1)$ is
already known by Theorem \ref{8rnotinnn}.) So now $t_k > -1$, and $v_t$
is maximal of highest weight in $Y(t_k)$. Then we have\\

\begin{prop}\label{P4}
$Y(t_k) = Z(t)$ (and $t \notin \nn$).
\end{prop}

\begin{proof}
The same sort of reasoning, using linear combinations of $F^i Y^j$, is
used here. We are looking at $V \subset Y(t_k) \subset Z(t_k)$. So let us
assume that $v_x = p(Y,F) v_{t_k} \in V$. Thanks to the $H$-action, we
may assume that $v_x$ is in a single weight space. Again, we know $v_t =
p_{t_k -t}(Y,F) v_{t_k}$, so we may say w.l.o.g. that $v_x = p'(Y,F)v_t +
F^l q(Y,F)v_{t_k} \in V$, by the Euclidean algorithm (considering all
these as polynomials in $Y$). Here, we can choose $q$ to be monic in $Y$,
and we of course have $l>0$ and deg$(q) < t_k - t$ (thereby splitting
$v_x$ into the ``higher degree" and ``$Z(t)$" components).\\

The key fact to be shown is that $q=0$. Suppose not, and let $v_x$ be a
vector in $V$ of highest weight $x$ for which $q \neq 0$. Now, we see
that $Ev_x = Ep' \cdot v_t + EF^l q \cdot v_{t_k} \in V$, and the second
term equals $(([E,F]F^{l-1} + \dots + F^{l-1}[E,F])q + F^l Eq) \cdot
v_{t_k} = F^{l-1}(\la + FE) q \cdot v_{t_k}$ for some scalar $\la$.
Clearly, this is in the \usl-span of the vectors $1, Y, \dots, Y^{t_k - t
- 1}$ (inside $Z(t_k)$) by Lemma \ref{Lbasis1}, since $q$ is monic. Hence
by maximality of weight of $v_x$, this second term is zero, because the
other term $Ep' \cdot v_t$ is in $H_f v_t$ (as $v_t$ is maximal).\\

Thus, $EF^l q \cdot v_{t_k} = 0$. But here, $l > 0$, so by Proposition
\ref{12general2} we know that $F^l q \cdot v_{t_k} = F^{j+1}v_j$ for some
$j$. Now look at $X^{-x-1}v_x \in V$. Since $t<-1$, hence the first term
of $v_x$ goes to $X^{-x-1}p' \cdot v_t \in Z(t)_{-1} = 0$. Thus
$X^{-x-1}v_x = X^{-x-1}F^{j+1}v_j$ and this has weight $-1$. Thus
$X^{-x-1}v_x = c_0 v_{-1}$ for some nonzero scalar $c_0$, so that $v_{-1}
\in V \subset Y(t_k)$. This is impossible, and hence $q = 0$ to start
with.
\end{proof}\hfill

Let us look at composition series now. We can directly see that $Y(r) /
Z(s-1) = Y(t_0) / Z(t_1) = Y(t_1,t_0) / Y(t_1,t_1)$ is simple (by
Proposition \ref{P3} above), and has highest weight vector $F^{s+1}v_s$.
Again, $Y F^{j+1}v_j = F^{j+1} (v_{j-1} - d_{j-1} Fv_{j+1})$, so we claim
inductively that $F^{j+1}v_j$ lies in $\mf{U}(N_-) (F^{s+1}v_s)$. This
holds in the base case because $v_{s-1} = 0$ in the simple quotient
$V(r,s)$.

Therefore $Y(r) / Z(s-1)$ is a simple standard cyclic module with highest
weight vector $F^{s+1}v_s$, hence of highest weight $-s-2$. So it is
isomorphic to $V(-s-2)$. We can now go to ``lower" $t_i$'s, and easily
calculate the composition factors.\\

Thus $Z(r)$ has a finite composition series. The set of composition
factors is $V(t_0), V(-t_1 -3), V(t_1), \dots, V(-t_k -3), V(t_k)$, and
the set of composition factors of $Y(t_k)$ (which is $0$ or $Z(t)$ from
above). If $Y(t_k) = Z(t)$ or $t_k = -1$ then we know everything about
the composition series of $Z(t_k)$, from Theorem \ref{8rnotinnn}. Thus,
in either case we know the composition factors of $Z(r)$ completely,
modulo the following remark.

\begin{remark}
The only question that needs answering is: Given $r,t_k$ as above, when
is $Z(t_k)$ simple ?\\
\end{remark}

If $r \notin \nn$ then there is only one Jordan-Holder series, and we
know {\it all} submodules of $Z(r)$. If $r \in \nn$, then there may be
more than one series; one example is
$$Z(r) = Y(t_0,t_0) \supset Y(t_1,t_0) \supset Y(t_1,t_1) \supset \dots
\supset Y(t_k,t_k) = Z(t_k) \supset Y(t_k) (\supset \dots)$$

where $Y(t_k) = Z(t)$ or 0. We thus have show the analogue of Theorem
\ref{8rnotinnn}, namely

\begin{theorem}\label{13rinnn}
Suppose $r \in \nn$, and $r = t_0 > t_1 > \dots > t_k \geq -1$ are the
various roots (in $\Z$) of $\aaa_{r,r-t+1}$.
\begin{enumerate}
\item The submodules of $Z(r)$ with highest weight vector of weight $\geq
-1$ are of the form $Y(t_i,t_j)$.
\item If $t_k > -1$, then either $Z(t_k)$ is simple, or $Y(t_k)$ has a
maximal vector of (highest) weight $t<-1$, whence $Y(t_k) = Z(t)$. In
this case, or if $t_k = -1$, we know the rest of the submodules from
Theorem \ref{8rnotinnn}.
\item $Z(r)$ has a finite composition series, of length at most
$4(\emph{deg}(f)+1)$.
\item The composition factors are simple modules $V(\la)$ with highest
weights $\{ t_i,\ -t_{i+1} - 3 : 0 \leq i \leq k-1 \}$ and $t_k$ if
$Z(t_k)$ is simple. If $Y(t_k) = Z(t)$, then we add the composition
factors of $Z(t)$ to this. Each simple module occurs with multiplicity 1
or 2.
\end{enumerate}
\end{theorem}

\noindent Thus, we can find all simple modules and primitive ideals in
this case. We can make similar claims for any $Z(r) \to V \to 0$ (where
$r \in \nn$). Some of the multiplicities may be 2, as we shall see
below.\\

\section{The (finite) sets $S(r)$ satisfy all the assumptions}

We are now ready to show that all the assumptions (and hence the
analysis) in the first part of the paper, hold in the case of $H_f$.

\begin{lemma}\label{finlength}
Every Verma module $Z(r)$ has finite length, so \calo = \caloo.
\end{lemma}

\begin{proof}
This follows from Lemmas \ref{caloo} and \ref{L2.2}.
\end{proof}

Thus, the assumptions and results of Theorem \ref{16finitelength} hold in
this case. Therefore every module in $\calo$ has an SC-filtration, is of
finite length, and $\calo$ is an abelian category that is self-dual as
well.

\begin{theorem}
If $Z(r)$ has a simple subquotient $V(t)$, then $S(r) = S(t)$.
\end{theorem}

\begin{proof}
This follows from Theorems \ref{8rnotinnn} and \ref{13rinnn}, since we
now explicitly know what composition factors any given Verma module can
have.
\end{proof}

\begin{remark}
Thus the $S(r)$'s decompose into a disjoint union of subsets, each of
which is finite, and plays the role of the $S(\la)$'s of the first part
of this paper. (We shall see below that in most cases the $S(r)$'s are
irreducible - and hence of the form $S(\la)$.)\\
\end{remark}

Over here, just as in the first part, we do not have the classical notion
of blocks. However, we can construct blocks as in the first part (using
the connected components of the $S(r)$'s), because all the assumptions
now hold. We define the block $\calo(r)$ to consist of all $M \in \calo$,
all of whose simple subquotients are of the form $V(t)$ for some $t \in
S(r)$.

Now all the results mentioned above hold, and we have enough projectives,
progenerators, and BGG reciprocity in the highest weight category
$\calo(r)$. We also have $\calo = \bigoplus \calo(r)$.\\

\section{More on the roots of $\aaa_{rt}$}

\noindent We actually know more about the roots of $\aaa_{rt}$, from the
following proposition.

\begin{prop}\label{P2.1}
\hfill
\begin{enumerate}
\item For all $r \in k,\ c_r = c_{-r-2}$, and hence $c_{0r} =
c_{0,-r-2}$.
\item $\aaa_{r,2r+4} = 0$ if $r+1 \in \nn$.
\item Suppose $r+1 \in \nn$. Then $Z(r)_{-2}$ has a maximal vector iff
$\aaa_{r,r+2} = 0$, iff $Z(r)_{-1}$ has a maximal vector.
\item If $r \notin \nn$ then the roots of $\aaa_{rt}$ in $r-\nn$ are
finitely many, as seen above.\\
If $r \in \nn$, then let $r_0$ be maximal in $S(r)$. Suppose $r_0 = t_0 >
\dots > t_k \geq -1$ are all roots of $\aaa_{r_0,r_0-t+1}$ in $r_0 - \nn
\cap \nn - 1$. Then the roots of $\aaa_{r,r-t+1}$ in $r-\nn$ are all
$t_j$'s less than $r$, and $\{ -t_j - 3 : 0 \leq j \leq k \}$. 
\item The length of any Verma module $Z(r)$ is at most $3\emph{deg}(f) +
4$.\\
\end{enumerate}
\end{prop}

\begin{remark}
If $Z(r)_{-1}$ has a maximal vector ($r \in \nn$) then $\aaa_{r,r+2} = 0$
and from part (3) above we see that $Z(r)_{-2}$ also has a maximal
vector. In this case, Corollary \ref{C9} seems to, but {\it does not}
imply, that $\mf{U}(N_-) v_{-2} \hookrightarrow \mf{U}(N_-)v_{-1}
\hookrightarrow Z(r)$. It may happen, actually, that $\mf{U}(N_-) v_{-2}
\subset Z(r) \supset \mf{U}(N_-)v_{-1}$, but $\mf{U}(N_-) v_{-2}
\nsubseteq \mf{U}(N_-)v_{-1}$. The reason this does not go through, is
that $d_{-3}$ is not defined.

Also note that not all multiplicities are zero; in particular, if $r_0$
is maximal in $S(r)$, then every single $V(t)$ (for $t \in S(r)$), except
at most for $V(-r_0 - 3)$, is a subquotient of $Z(r_0)$. Further, part
(5) holds for any $Z(r) \to V \to 0$, and is a better estimate than
above.

Next, we observe that if a block $S(r) \subset \Z$ has size 2, then it
may not be irreducible, as in the original definition of $S(\la)$ (in the
general case) ! In this case, we work with each element as a block by
itself. But in all other cases, each set $S(r)$ is a block by itself
(i.e. ``irreducible", as in the first part). This follows from the
remarks above, and Theorems \ref{8rnotinnn} and \ref{13rinnn}.

Finally, observe that if $D$ is the unipotent decomposition matrix, then
each entry of $D$ is 0, 1 or 2, as we saw in the separate cases $r \in
\nn$ and $r \notin \nn$ above.\\
\end{remark}

\begin{proof}
(1), (2) and (4) are calculations. As for (3), one way is clear, by
Theorem \ref{Tnec}. Conversely, suppose $\aaa_{r,r+2} = 0$. Then we can
verify that $v_{-2} = Yv_{-1} - c_{0,-1}Fv_0$ is indeed a maximal
vector.\\

\noindent (5) For $r \notin \nn$ this is clear from Theorem
\ref{8rnotinnn}. For $r \in \nn$ we recall the structure of $Z(r)$. We
know from the previous part, that $n_+ \geq k$. Here, we define $n_+$ to
be the number of roots of $\aaa_{rt}$ (out of a total of $2k+2$ roots, as
given), that are in $\nn$.

Thus the number of negative integer roots $n_-$ is at most $k + 2$. There
are at most two simple subquotients (in $Y(-3-t)$ and then in $Z(t)$, as
earlier) for each of these, and one simple subquotient for each positive
root.

Hence the total number of terms in a composition series is at most $2n_-
+ n_+ = (n_- + n_+) + n_- \leq (2k+2) + (k+2) = 3k+4$. But $2k+2 \leq 2
\deg(f)+2$ by Corollary \ref{C11}, so $k \leq$ deg$(f)$, whence the
length of a composition series is $\leq 3k+4 \leq 3\deg(f) + 4$, as
claimed.
\end{proof}\hfill

\begin{remark}
It remains to find out the composition series of a Verma module for the
case $r \in \nn$, or equivalently, the composition series for $Z(t_k)$ in
this case. This would lead to a complete knowledge of all multiplicities
$[Z(\la) : V(\mu)]$. However, we do not know the answer to this
question.

One guess would be that $Z(t) \hookrightarrow Z(r)$ iff $\aaa_{r,r-t+1} =
0$, since one implication holds in general, and the other holds as well,
if $r \notin \nn$. However, this converse implication is {\it false} for
$r \in \nn$. For example, setting $g(T) = 1+f(T)$, direct calculations
yield that when $t_k = -1,\ Z(-2) \hookrightarrow Z(-1)$ iff $c_{0,-1} =
g(-1/8) = 0$. Similarly, when $t_k = 0,\ Z(-3) \hookrightarrow Z(0)$ iff
$g(0)(g'(0) / 2 + g(-1/8)) = 0$, and this is not true for general $g$
(e.g. $g=1$, or $f=0$).\\
\end{remark}

\section{Weyl's theorem fails, multiplicities may be 2, and more}

We now look at a specific module $Z(0)$. Suppose $f$ has the property
that $c_{00} = c_{0,-1} = 0$. Then $Z(0)$ has maximal vectors $v_0,
v_{-1}, v_{-2}, v_{-3}$, and $v_i = Y^i v_0$ for each of these.\\

Observe that in general, we cannot obtain a resolution for $V = V(r,s)$
in terms of the $Z(\la)$'s. In any such resolution, the first term would
be $Z(r) \twoheadrightarrow V(r,s)$. We then need some $\mu$ so that
$Z(\mu) \twoheadrightarrow Y(r)$. But this is not true in general: look
at the above example $V = Z(0)$. Clearly, $Z(0) \twoheadrightarrow
V(0,0)$ has kernel $Y(0) = (Y,F)$. Clearly, if $\varphi : Z(\mu)
\twoheadrightarrow Y(0)$, then $v_\mu \mapsto Y$ (for if it maps to zero,
then $\varphi = 0$). But then we see that $F \notin \image(\varphi)$.\\

Also, observe that the multiplicities $[Z(r) : V(r')]$ are not 0 or $1$
in general: in the above example, we see that $[Z(0) : V(-2)] = 2$. This
is because we have the series $Z(0) \supset Y(0) = (F,Y) \supset Z(-1) =
(Y) \supset Y(-1) = Z(-2) = (Y^2) \supset Y(-2) = Z(-3) = (Y^3) \supset
Y(-3) \supset \dots$, and the subquotients are $V(0), V(-2), V(-1),
V(-2), V(-3), \dots$.\\

Finally, we provide a counterexample to Weyl's theorem - namely, a
(finite-dimensional) $H_f$-module $M$ and a submodule $N$ in it that has
no complement. Take $M = V(1,0) \supset V(0,0) = N$, i.e. $M = H_f / I$,
where the left ideal $I$ is generated by $(H-1),E,X,Y^2,FY,F^2$. In other
words, $M = k w_1 \oplus k w_0 \oplus k w_{-1}$, and $N = k w_0$, with
module relations as follows:
\begin{eqnarray*}
&& E w_1 = X w_1 = 0; \qquad F w_{-1} = Y w_{-1} = 0;\\
&& F w_1 = w_{-1},\  E w_{-1} = w_1; \qquad Y w_1 = w_0,\ X w_{-1} = -w_0
\end{eqnarray*}

and $X w_0 = Y w_0 = H w_0 = E w_0 = F w_0 = 0$ (i.e. $w_0$ is killed by
$X,Y,E,F,H$).

It can be checked that this is a valid $H_f$-module structure on $M$, if
we have $c_{00} = c_{01} = 0$. However, it is obvious that $k w_0$ is a
submodule (with a trivial module structure). Any complement must contain
$w_1$ + lower weight vectors, but when we apply $Y$ to this, we get
$w_0$. Thus $w_0$ lies in the submodule and in its complement; a
contradiction. Hence there does not exist a complement to $kw_0$ in $M$,
and Weyl's theorem fails for this case.\\

\appendix
\section{Algebraic preliminaries}

Throughout, $R$ denotes a ring, and \calo denotes an abelian subcategory
of $R$-mod.\\

\begin{prop}\label{extadditive}
If $0 \to A \oplus B' \to C \to B'' \to 0$ in \calo, and $\Ext^1(B'',A) =
0$, then $C = A \oplus B$, where $0 \to B' \to B \to B'' \to 0$ in
\calo.
\end{prop}

\begin{proof} Apply $\hhom(B'',-)$ to the s.e.s. $0 \to B' \to B' \oplus
A \to A \to 0$. Then our result follows by considering the long exact
sequence of \Ext's.
\end{proof}\hfill

\begin{prop}\label{extvecspc}
Suppose $R$ is a $k$-algebra, where $k$ is a field, and say we have an
exact contravariant duality functor $F : \calo \to \calo$ (i.e. $F(M)
\subset \hhh_k(M,k),\ F(F(M)) = M$). Then $F : \Ext^1(M'',M') \to
\Ext^1(F(M'),F(M''))$ is an isomorphism of $k$-vector spaces.
\end{prop}

The proof more or less follows from the way we define the vector space
operations; they use pullbacks, push-forwards, and element chasing in
commutative diagrams, e.g. cf. \cite{F}.\\

\noindent \textbf{Setup :} Now suppose also that \calo is finite length,
and a full subcategory of $R$-mod.
Let \calp denote all indecomposable projective objects in \calo, and
let \cals denote all simple objects. (Thus Fitting's Lemma holds.)

\begin{theorem}\label{fitting}
\hfill
\begin{enumerate}
\item Every object $P$ in \calp has a unique maximal sub-object
$(\Rad(P)).\ P$ is the projective cover of $P / \Rad(P) \in \cals$.
\item The map $F : \calp \to \cals$ given by $F(P) = P / \Rad(P)$ is
one-one. If enough projectives exist in \calo, then $F$ is a bijection.
\end{enumerate}
\end{theorem}

\begin{theorem}\label{bass}
Suppose now that enough projectives exist in \calo, and \calp is finite.
\begin{enumerate}
\item $Q = \bigoplus_{P \in \calp} P^{\oplus n_P}$ is a progenerator for
\calo, as long as all $n_P \in \N$.
\item Set $B = \hhom(Q,Q)$. Then $B$ is unique upto Morita equivalence,
and the functor $D = \hhom(Q,-)$ is an equivalence between $\calo$ and
$(\emph{mod-}B)^{fg}$ (i.e. finitely generated right $B$-modules).
\item $D$ and $E = Q \otimes_B -$ are inverse equivalences between \calo
and $(\emph{mod-}B)^{fg}$.
\end{enumerate}
\end{theorem}

\noindent (Part (2) of Theorem \ref{bass} is from \cite[Pg.
55]{Bass}.)\\

\vspace{3ex}

\noindent {\bf Acknowledgements :} I thank my advisor Prof. Victor
Ginzburg for suggesting the problem to me, as well as for his help and
guidance.\\


\providecommand{\bysame}{\leavevmode\hbox to3em{\hrulefill}\thinspace}
\providecommand{\MR}{\relax\ifhmode\unskip\space\fi MR }
\providecommand{\MRhref}[2]{%
  \href{http://www.ams.org/mathscinet-getitem?mr=#1}{#2}
}
\providecommand{\href}[2]{#2}

\end{document}